\documentclass[12pt]{statsoc}

\usepackage{amsmath} 
\usepackage{amstext}
\usepackage{amssymb}
\usepackage{enumerate}
\usepackage{slashbox}
\usepackage{array}
\usepackage{color}
\usepackage{bbm}
\usepackage{natbib}
\usepackage{graphicx}
\usepackage{caption}
\usepackage[title]{appendix}
\usepackage{url} % not crucial - just used below for the URL 
\usepackage{fancyhdr}
\usepackage{float}

\def\bs{\boldsymbol}

\newtheorem{thm}{Theorem}[section]
\newtheorem{assumption}{Assumption}
\newtheorem{rem}{Remark}

\newcommand{\argmax}{\operatornamewithlimits{argmax}}

\def\boxit#1{\vbox{\hrule\hbox{\vrule\kern6pt
			\vbox{\kern6pt#1\kern6pt}\kern6pt\vrule}\hrule}}

\title{Nonparametric Estimation of the Continuous Treatment Effect with Measurement Error} 
\author{Wei Huang}
\address{School of Mathematics and Statistics,
University of Melbourne, 
Australia.}\email{wei.huang@unimelb.edu.au}
\author[W. Huang and Z. Zhang]{Zheng Zhang} 
\address{Center for Applied Statistics, Institute of Statistics \& Big Data, 
Renmin University of China,
China.}\email{zhengzhang@ruc.edu.cn}

\allowdisplaybreaks
\begin{document}
\sloppy

\begin{abstract}
We identify the average dose-response function (ADRF) for a continuously valued error contaminated treatment by a weighted conditional expectation. We then estimate the weights nonparametrically by maximising a local generalised empirical likelihood subject to an expanding set of conditional moment equations incorporated into the deconvolution kernels. Thereafter, we construct a deconvolution kernel estimator of ADRF. We derive the asymptotic bias and variance of our ADRF estimator and provide its asymptotic linear expansion, which helps conduct statistical inference. To select our smoothing parameters, we adopt the simulation-extrapolation method and propose a new extrapolation procedure to stabilise the computation. Monte Carlo simulations and a real data study illustrate our method's practical performance.
\end{abstract}

\vfill

\footnotetext{\footnotesize{\textit{The authors contributed equally to this work and are listed in the alphabetical order.}}}
\newpage

\section{Introduction}
Identifying and estimating the causal effect of a treatment or policy from observational studies is of great interest to economics, social science, and public health researchers. There, confounding issues usually exist (i.e. individual characteristics are related to both the treatment selection and the potential outcome), making the causal effect not directly identifiable from the data. 
Early studies focused on whether an individual receives the treatment or not (e.g. \citealp{Rosenbaum1983,rosenbaum1984reducing,hahn1998role,Hirano2003}). More recently, as opposed to such binary treatments, researchers have been investigating the causal effect of a continuously valued treatment, where the effect depends not only on the introduction of the treatment but also on the intensity; see \cite{hirano2004propensity, galvao2015uniformly, kennedy2017non, Fong_Hazlett_Imai_2018, huber2020direct, dong2019regression, Ai_Linton_Motegi_Zhang_cts_treat, Ai2020}, and \cite{huang2021unified}, among others. However, all these methods require the treatment data to be measured without errors.

In many empirical applications, the treatment and confounding variables may be inaccurately observed. For example, \cite{mahajan2006identification} studied the binary treatment effect when the observed treatment indicators are subject to misclassification; see also \cite{lewbel2007estimation} and \cite{molinari2008partial}, among others. \cite{battistin2014treatment} investigated the extent to which confounder measurement error affects the analysis of the treatment effect. Continuous treatment variables are also likely to be measured with error in practice. For example, in the Epidemiologic Study Cohort data from the first National Health and Nutrition Examination (NHANES-I) (see \citealp{Carroll2006}), over 75\% of the variance in the fat intake data is made up of measurement error.  However, to the best of our knowledge, no existing work has considered the identification and estimation of the causal effect from such error-contaminated continuous treatment data. To bridge this gap in the body of knowledge, we focus on continuous treatment data measured with classical error; that is, instead of observing the received treatment, researchers only observe the sum of the treatment and a random error. 

In particular, we study a causal parameter of primary interest to scholars in the literature, the average dose-response function (ADRF). It is defined as the population mean of the individual potential outcome corresponding to certain levels of the treatment.
To resolve the confounding problem, we adopt the most widely imposed assumption in the literature on treatment effects, the unconfoundedness condition (e.g. \citealp{Rosenbaum1983, rosenbaum1984reducing, hirano2004propensity}), which assumes that the assignment of the treatment is independent of the potential outcome of interest given a set of observable covariates. 

Perhaps, the most straightforward approach to this problem is based on some parametric model specifying how the outcome relates to confounders and the treatment. Then one may apply the parametric methods to measurement error data in the literature (see~e.g.~\citealp{Carroll2006}). However, the parametric approach suffers from the model misspecification problem and does not incorporate available information on the treatment mechanism. Thus, this paper focuses on robust nonparametric estimation of ADRF.

In the literature on error-free continuous treatment data, nonparametric estimation of ADRF under unconfoundedness assumption has been extensively studied. For example, \cite{galvao2015uniformly} identified the ADRF via an unconditional weighted expectation, with the weighting function being the ratio of two conditional densities of the treatment variable. They estimated these two conditional densities separately and then constructed the estimator for the ADRF. However, it is well known that such a ratio estimator can be unstable owing to its sensitivity to the denominator estimation \citep{kang2007demystifying}. To improve the robustness of estimation, \cite{kennedy2017non} developed a doubly robust estimator for ADRF by regressing a doubly robust mapping on the treatment; see more detailed discussion in Remark~2 of Section~\ref{sec:estiamte_pi}. \cite{Ai_Linton_Motegi_Zhang_cts_treat} identified ADRF via a weighted conditional expectation on the treatment variable, where the weighting function is the ratio of the treatment variable's marginal density and its conditional density given the confounders. They estimated the weighting function (but not the two densities in the ratio separately) by maximising  entropy subject to some moment restrictions  that identify the weighting function, then obtained the nonparametric weighted estimator of the ADRF.  The idea of  estimating the density ratio directly has also been exploited in the literature in (bio)statistics based on parametric modeling \citep{qin1998inferences}, semiparametric modeling \citep{cheng2004semiparametric}, and machine learning with an augmented dataset  \citep{diaz2021nonparametric}.

However, no existing methods for the error-free treatment data can be easily extended to the error-contaminated treatment data. Indeed, under the unconfoundedness assumption, all those methods depend on estimating a weighting function with the treatment's conditional density given  the confounders in the denominator. Its estimation is a critical step and full of challenges in the presence of measurement error. Although the nonparametric method of estimating the (conditional) density of a variable measured with the classical error has been developed (see e.g.~\citealp{ fan1991asymptotic, fan1991optimal, meister2006density}), estimating the density and plugging it into the denominator suffers from the instability issue mentioned above. Such an issue can be even worse when the treatment is measured with errors. Unfortunately, existing techniques of improving the robustness cannot apply to the measurement error data; for example, \citeauthor{Ai_Linton_Motegi_Zhang_cts_treat}'s \citeyear{Ai_Linton_Motegi_Zhang_cts_treat} approach heavily relies on the moment information of a general function of the treatment variable, which is hard to obtain in the presence of measurement error.\label{ReplyContribution1} To our best knowledge, estimating the moment of a general function of a variable measured with errors and its associated theoretical behaviour are challenging problems and remain open in the literature (see \citealp{hall2008estimation} and our discussion in section~\ref{sec:estiamte_pi} for more details).

We propose a broad class of novel and robust nonparametric estimators for the ADRF when the treatment data are measured with error. We first represent the ADRF as in \cite{Ai_Linton_Motegi_Zhang_cts_treat}, i.e.~a weighted expectation of the observed outcome on the treatment variable, where the weighting function is a ratio of the treatment variable's marginal density and its conditional density given the confounders. Then we propose a novel approach to identify the weighting function with no need to estimate the moment of a general function of the treatment. Specifically, we propose a \emph{local} empirical likelihood technique that only requires estimating the conditional moment of a general function of the error-free confounders given each fixed value of the treatment. A consistent nonparametric estimator of such a conditional moment can be obtained by using a local constant deconvolution kernel estimator in the literature on errors-in-variables regression; see \cite{Fan1993, carroll2004low, hall2007ridge, delaigle2009design}, and \cite{Meister2009}, among others. Moreover, by choosing different criterion functions for the local generalised empirical likelihood, the proposed method produces a wide class of weights, including exponential tilting, empirical likelihood and generalised regression as important special cases. Once the weights are consistently estimated, we construct our estimator of the ADRF by estimating the weighted conditional expectation of the observed outcome using another local constant deconvolution kernel estimator.\label{ReplyContribution2}

Based on our identification theorem and the literature on the local constant deconvolution kernel estimator, it is not hard to see the consistency of our estimator. However, the asymptotic behaviour of a local constant deconvolution kernel estimator incorporated with such local generalised empirical likelihood weight has never been investigated. %The theories for the global generalised empirical likelihood framework do not apply.
We show the asymptotic bias of our estimator consists of two parts: the one from estimating the weights and the one from estimating the weighted conditional expectation, but it can achieve the same convergence rate as that in the error-free case. The asymptotic variance depends on the type of measurement error. We study the most commonly used two types of measurement errors: the ordinary smooth and the supersmooth ones. The convergence rate of the variance is slower than that in the error-free case but can achieve the optimal rate for the nonparametric estimation in measurement error literature. 

Moreover, without imposing additional conditions, it is well known that the explicit expression of the asymptotic variance of the deconvolution kernel for supersmooth error (e.g.~Gaussian error) is not derivable; even the exact convergence rate remains unknown. Therefore, it is common to use residual bootstrap methods to make the inference (see e.g.~\citealp{fan1991asymptotic, delaigle2009design,Delaigle2015}). However, residual bootstrap methods  do not work well in our treatment effect framework. We thus provide the asymptotic linear expansion of our ADRF estimator. Based on this, we propose an undersmoothing (pointwise) confidence interval for the ADRF.  

A data-driven method to select the bandwidth for a local polynomial deconvolution kernel estimator is challenging. The only method in the literature to our best knowledge is the simulation-extrapolation (SIMEX) method proposed by \cite{Delaigle2008SIMEX}. 
However, the linear back-extrapolation method suggested by the authors performed unstably in our numerical studies. We, thus, propose a new local constant extrapolant function. 
Monte Carlo simulations show that our estimator is stable and outperforms the naive estimator that does not adjust to the measurement error. We also demonstrate the practical value of our estimator by studying the causal effect of fat intake on breast cancer using the Epidemiologic Study Cohort data from NHANES-I (\citealp{Carroll2006}).

The remainder of the paper is organised as follows. We introduce the basic framework and notations in Section~\ref{sec:DataandModel}. Section~\ref{sec:estimation} presents our estimation procedure, followed by the asymptotic studies of our estimator in Section~\ref{sec:asymptotic}. We discuss the method of selecting our smoothing parameters in Section~\ref{sec:ParameterSelection}. Finally, we present our simulation results and real data analysis in Section~\ref{sec:Simulation}.

\section{Basic Framework}\label{sec:DataandModel}
We consider a continuously valued treatment in which the observed treatment variable is denoted by $T$ with the probability density function $f_T(t)$ and support $\mathcal{T}\subset \mathbb{R}$. Let $Y^\ast(t)$ denote the potential outcome if one was treated at level $t$ for $t\in\mathcal{T}$. In practice, each individual can only receive one treatment level $T$ and we only observe the corresponding outcome $Y:=Y^\ast(T)$. We are also given a vector of covariates $\boldsymbol{X}\in\mathbb{R}^r$, with $r$ a positive integer, which is related to both $T$ and $Y^\ast(t)$ for $t\in\mathcal{T}$. Moreover, we consider the situation in which the treatment level is measured with classical error, that is, instead of observing $T$, we observe $S$ such that
\begin{equation}
S = T + U,\label{ME}
\end{equation}
where $U$ is the measurement error, independent of $T, \bs{X}$ and $\{Y^\ast(t)\}_{t\in\mathcal{T}}$, and its characteristic function $\phi_U$ is known; see Remark \ref{rem:UnknownfU} for the case in which $\phi_U$ is unknown.

The goal of this paper is to nonparametrically estimate of the unconditional ADRF $\mu(t):= \mathbb{E}\{Y^\ast(t)\}$ from an independent and identically distributed (i.i.d.) sample $\{S_i,\bs{X}_i,Y_i\}_{i=1}^N$ drawn from the joint distribution of $(S,\bs{X},Y)$. %Then methods of testing distributional effects, i.e. $H_0: F_{t_0}(y)=F_{t_1}(y)$ for all $y\in\mathcal{Y}$ v.s. $H_1: F_{t_0}(y)\neq F_{t_1}(y)$ for some $y\in\mathcal{Y}$, for any two treatment levels $t_0\neq t_1$, will be established.

%This is a test of common interest in causal inference literature. It is testing the difference between the potential outcomes of applying treatment $t_0$ and $t_1$ to a same individual. Note that this is different to testing if $P(Y\leq y|T=t_1) = P(Y\leq y|T=t_0)$ for all $y\in \mathcal{Y}$, which essentially tests the effect difference between the subgroup of people who have received treatment $t_1$ and those who have received $t_0$.

$Y^\ast(t)$ is never observed simultaneously for all $t\in\mathcal{T}$ or even on a dense subset of $\mathcal{T}$ for any individual, but only at a particular level of treatment, $T$. Thus, to identify $\mu(t)$ from the observed data, the following  assumption is imposed in most of the treatment effect literature (e.g. \citealp{rubin1990comment, kennedy2017non, Ai_Linton_Motegi_Zhang_cts_treat,Ai2020, d2021overlap}). 
\begin{assumption} We assume \label{UnconfoundAssump}
\begin{enumerate}
	\item [(i)]  	(Unconfoundedness)
	for all $t\in\mathcal{T}$, given $\boldsymbol{X}$, $T$ is independent of $Y^*(t)$, that is, $T \perp \{Y^\ast(t)\}_{t\in\mathcal{T}} | \bs{X}$;
	\item [(ii)] %{\color{red}(Stable Unit Treatment Value Assumption) Let $\boldsymbol{T}=(T_1,\ldots, T_N)$ be the vector of treatment assignments to individuals $i=1,\ldots,N$. For any possible assignments $\boldsymbol{T}$ and $\boldsymbol{T}'$, if $T_i=T_i'$, then $Y_i^\ast(T_i) = Y_i^\ast(T_i')$, for $i=1\ldots,N$; }
	(No Interference) For $i=1,\ldots, N$, the outcome of individual $i$ is not affected by the treatment assignments to any other individuals. That is, $Y_i^\ast(T_i,\boldsymbol{T}_{(-i)})=Y_i^\ast(T_i,\boldsymbol{T}'_{(-i)})$ for any $\boldsymbol{T}_{(-i)}$, $\boldsymbol{T}'_{(-i)}$, where $Y_i^\ast(T_i,\boldsymbol{T}_{(-i)})$ is the potential outcome of individual $i$ given the treatment assignments to individual $i$ and the others are $T_i$ and $\boldsymbol{T}_{(-i)}:=(T_1,\ldots, T_{i-1}, T_{i+1},\ldots, T_N)$, respectively;
%	\item [(ii)] {\color{red}(Stable Unit Treatment Value Assumption) $T_i =t$ implies $Y_i = Y^\ast_i(t)$ for $i=1,\ldots, N$;}
	\item [(iii)] (Consistency) $Y=Y^*(t)$ a.s. if $T=t$;
	\item [(iv)] (Positivity) the generalised propensity score $f_{T|X}$ satisfies $f_{T|X}(t|\boldsymbol{X})>0$ a.s. for all $t\in\mathcal{T}$.
\end{enumerate}
\end{assumption}
%{\color{red}Assumption~\ref{UnconfoundAssump} (ii) \citep{rubin1990comment} comprises two subassumptions \citep{rubin2005causal}: 1) the outcome of each individual $Y^\ast_i(t)$ is not affected by what assignment any other individual received (i.e.~no interference); 2) there is no hidden versions of treatments so that no matter how individual $i$ received treatment $T_i$, the observed outcome would be $Y_i^\ast(T_i)$ (i.e.~consistency).} 
Under Assumption \ref{UnconfoundAssump}, for every fixed $t\in\mathcal{T}$,  $\mu(t)$ can be identified as follows:
\begin{align}
	\mu(t) 
	= &\mathbb{E}[\mathbb{E}\{Y^\ast(t)|\bs{X}\}]
	= \mathbb{E}[\mathbb{E}\{Y^\ast(t)|\bs{X}, T=t\}] \quad \text{(by Assumption \ref{UnconfoundAssump}  (i) and (ii))}\notag \\
	= &\mathbb{E}\{\mathbb{E}(Y|\bs{X},T=t)\} \quad \text{(by Assumption \ref{UnconfoundAssump}  (iii))}\notag\\
    =& \int_{\mathcal{X}}\int_{\mathcal{Y}}\frac{f_{T}(t)}{f_{T|X}(t|\bs{x})}yf_{Y|X,T}(y|\bs{x},t)f_{X|T}(\bs{x}|t)\,dy\,d\bs{x}  \quad \text{(by Assumption \ref{UnconfoundAssump}  (iv))} \nonumber\\
%=& \int_{\mathcal{X}}\int_{\mathcal{Y}}\frac{f_{T}(t)}{f_{T|X}(t|\bs{x})}yf_{Y,X|T}(y,\bs{x}|t)\,dy\,d\bs{x}\nonumber\\
	=& \mathbb{E}\{\pi_0(t,\bs{X})Y|T=t\}\,,\label{PC}
\end{align}
where 
\begin{equation}
	\pi_0(t,\bs{x}) := \frac{f_T(t)}{f_{T|X}(t|\bs{x})}\,,\label{pi0def}
\end{equation}
with $f_T$ and $f_{T|X}$ being the density function of $T$ and the conditional density of $T$ given $\bs{X}$, respectively. The function $\pi_0(t,\boldsymbol{x})$ is called the \emph{stabilised weights} in \cite{Ai_Linton_Motegi_Zhang_cts_treat}. 

If $T$ is fully observable and $\pi_0(t,\boldsymbol{x})$ is known, estimating $\mu(t)$ in \eqref{PC} is reduced to a standard regression problem. For example, $\mu(t)$ can be consistently estimated using the Nadaraya--Watson estimator:
\begin{equation}
	\mu_{NW}(t):=\frac{\sum^N_{i=1}\pi_0(t,\boldsymbol{X}_i)Y_iL\{(t-T_i)/h\}}{\sum^N_{i=1}L\{(t-T_i)/h\}}, \quad t\in\mathcal{T}\,,\label{LocalConstant}
\end{equation}
where $L(\cdot)$ is a prespecified univariate kernel function such that $\int_{-\infty}^\infty L(x)dx=1$, $h$ is the bandwidth. However, we do not observe $T$ but only $S$ in \eqref{ME} and $\pi_0(t,\bs{x})$ is also unknown in practice. We address these issues in the next section. 

\begin{rem}\label{rem:UnknownfU}
	The characteristic function of the measurement error $\phi_U$ may be unknown in some empirical applications. Several methods of consistently estimating it have been proposed in the literature. For example, \cite{diggle1993fourier}  assumed that the data from the error distribution are observable and proposed estimating $\phi_U$ from the error data nonparametrically. In some applications, the error-contaminated observations are replicated, and we can estimate the density from these replicates; see \cite{delaigle2009design}. When $\phi_U$ is known up to certain parameters, parametric methods are applicable (e.g. \citealp{meister2006density}). Finally, a nonparametric method without using any additional data is also available (e.g. \citealp{delaigle2016methodology}). Once a consistent estimator of $\phi_U$ is obtained, our proposed method can be directly applied.
	
	In particular, the assumptions that, %(a) data from the error distribution are available from other studies; 
	(a) the error distribution is known up to certain parameters that are identifiable from some previous studies and (b) repeated error-contaminated data are available, are commonly met in practice. 
	%If assumption~(a) is satisfied, we can estimate {\color{red}$\phi_U$} by $\widehat{\phi}_U(t) = N_U^{-1}\sum^{N_U}_{j=1}\exp(itU_j) $, {\color{red}which is a $\sqrt{N_U}$-consistent estimator with $N_U$} the sample size of the error data,  (\citealp{diggle1993fourier}). 
	When assumption~(a) is satisfied, we can usually obtain $\sqrt{N_U}$-consistent estimators of the parameters and thus a $\sqrt{N_U}$-consistent estimators of $\phi_U(t)$ uniformly in $t\in\mathbb{R}$, where $N_U$ is the sample size of the previous studies. For example, in our real data example in Section~\ref{sec:RealData}, the error distribution is known up to the variance. Since the convergence rate of our proposed estimator $\widehat{\mu}(t)$ is slower than $N^{-1/2}$, the asymptotic behaviour is insensitive to the estimation of $\phi_U$ provided that  $N=O(N_U)$. 
	Under assumption~(b), we observe 
	$$
	S_{jk} = T_{j} + U_{jk}, \quad k=1,2, \quad j=1,\ldots, N\,,
	$$
	where the $U_{jk}$'s are i.i.d.. Note that 
	$$\mathbb{E}[\exp\{it(S_{j1}-S_{j2})\}] = \mathbb{E}\{\exp(itU_{j1})\}\mathbb{E}\{\exp(-itU_{j2})\} = |\phi_U(t)|^2\,.
	$$ \cite{delaigle2008deconvolution} proposed to estimate $\phi_U(t)$ by $\widehat{\phi}_U(t) = \big|N^{-1}\sum^N_{j=1}\cos\{it(S_{j1}-S_{j2})\}\big|^{1/2} $. They showed that the deconvolution kernel local constant estimator using this $\widehat{\phi}_U(t)$ and that using $\phi_U(t)$ asymptotically behave the same under some regularity conditions. Specifically, when $U$ is ordinary smooth (as defined in \eqref{OrdinaryS}), they require $f_T$ to be sufficiently smooth relative to $U$'s density. In the contrast, if $U$ is supersmooth as in \eqref{SuperS}, the optimal convergence rate of a deconvolution kernel local constant estimator is logarithmic in $N$, which is so slow that the error incurred by estimating $\phi_U$ is negligible. This result can be extended to our setting. 
\end{rem}

\section{Estimation Procedure}\label{sec:estimation}
To overcome the problem that the $L\{(t-T_i)/h\}$'s are not empirically accessible, we apply the \emph{deconvolution} kernel approach (e.g. \citealp{Stefanski1990, Fan1993}). This method is often used in nonparametric regression in which the covariates are measured with classical error, as in \eqref{ME}, and the idea is introduced as follows. The density of $S$ is the convolution of the densities of $T$ and $U$, meaning that $\phi_S(w) = \phi_T(w)\phi_U(w)$, where $\phi_S$ and $\phi_T$ are the characteristic functions of $S$ and $T$, respectively. We consider $U$ with $\phi_U(w)\neq 0$ for all $w\in \mathbb{R}$. Using the Fourier inversion theorem, if $|\phi_T|$ is integrable, we have
\begin{align}\label{eq:fT_InverseF}
f_T(t) = \frac{1}{2\pi}\int_{-\infty}^\infty \exp(-iwt) \frac{\phi_S(w)}{\phi_U(w)}dw.
\end{align}
This inspired \cite{Stefanski1990} to estimate $f_T$ by $\widehat{f}_{T,h}(t):=(Nh)^{-1}\sum^N_{i=1}L_U\{(t-S_i)/h\}$, where
\begin{equation}
	L_U(v):= \frac{1}{2\pi} \int_{-\infty}^{\infty} \exp(-iwv)\frac{\phi_L(w)}{\phi_U(w/h)}dw,\label{KUdef}
\end{equation}
with $\phi_L$ the Fourier transform of the kernel $L$, which aims to prevent $\widehat{f}_{T,h}$ from becoming unreliably large in its tails.

Based on this idea, \cite{Fan1993} proposed a consistent errors-in-variables regression estimator by replacing the $L\{(t-T_i)/h\}$'s in \eqref{LocalConstant} with the $L_U\{(t-S_i)/h\}$'s. In our context, an errors-in-variables estimator of $\mu(t)$ is
\begin{equation}
\widetilde{\mu}(t) := \frac{\sum^N_{i=1}\pi_0(t,\bs{X}_i)Y_iL_U\{(t-S_i)/h\}}{\sum^N_{i=1}L_U\{(t-S_i)/h\}}.\label{OracleNonNormalise}
\end{equation}
Note that $U\perp (T,\boldsymbol{X},Y)$, we have
\begin{align}
	&\mathbb{E}\big[L_U\{(t-S)/h\}|T,\boldsymbol{X},Y \big]=\mathbb{E}\big[L_U\{(t-S)/h\}|T \big]\notag\\ 
	%&=& \mathbb{E}\bigg[\frac{1}{2\pi} \int_{-\infty}^{+\infty} e^{-iw(t-S)/h}\frac{\phi_L(w)}{\phi_U(w/h)}\,dw \bigg|T \bigg]\notag\\
	%&=&\frac{1}{2\pi} \int_{-\infty}^{+\infty} \mathbb{E}[\exp\{-iw(t-T-U)/h\}|T]\cdot\frac{\phi_L(w)}{\phi_U(w/h)}\,dw\notag\\
	=&\frac{1}{2\pi} \int_{-\infty}^{\infty}\exp\{-iw(t-T)/h\} \mathbb{E}[\exp(iwU/h)]\frac{\phi_L(w)}{\phi_U(w/h)}\,dw\notag\\
	=&\frac{1}{2\pi} \int_{-\infty}^{\infty} \exp\{-iw(t-T)/h\}\phi_L(w)\,dw%\notag\\
	=L\{(t-T)/h\},\label{LUCondExp}
\end{align}
where the last equation comes from the Fourier inversion theorem. Using this property, $\widetilde{\mu}(t)$ has the same asymptotic bias as that of $\mu_{NW}(t)$, which shrinks to zero as $h\to 0$. Then, to verify its consistency to $\mu(t)$, it suffices to show that its asymptotic variance decays to zero as $N\rightarrow\infty$, using a straightforward extension of the proof in \cite{Fan1993}.

However, the challenge is that $\pi_0$ is unknown in practice. We next show how to estimate $\pi_0(t,\bs{X})$ from the error-contaminated data $(S_i,\boldsymbol{X}_i,Y_i)$, $i=1,\ldots,N$.

\subsection{Estimating $\pi_0(t,\bs{X})$}\label{sec:estiamte_pi}
Observing \eqref{pi0def}, a straightforward way to estimate $\pi_0$ is to estimate $f_T$ and $f_{T|X}$ and then compute the ratio. However, this ratio estimator is sensitive to low values of $f_{T|X}$ since small errors in estimating $f_{T|X}$ lead to large errors in the ratio estimator (see~\citealp{Ai_Linton_Motegi_Zhang_cts_treat, Ai2020} for an example and Appendix~A.1.1 in the supplementary file for a detailed illustration).  As in the literature of error-free treatment effect, we treat $\pi_0$ as a whole and estimated directly to mitigate this problem. In paticular, we estimated it nonparametrically from an expanding set of equations, which is closely related to the idea in \cite{Ai_Linton_Motegi_Zhang_cts_treat}. However, their method is not applicable to error-contaminated data. 

Specifically, when the $T_i$'s are fully observable, \cite{Ai_Linton_Motegi_Zhang_cts_treat} found that the moment equation 
\begin{equation}
	\mathbb{E}\{\pi_0(T,\bs{X})u(\bs{X})v(T)\} = \mathbb{E}\{v(T)\}\mathbb{E}\{u(\bs{X})\}\label{UCM}
\end{equation}
holds for any integrable function $u(\bs{X})$ and $v(T)$, and that it identifies $\pi_0(\cdot,\cdot)$. They further estimated the function $\pi_0(\cdot,\cdot):\mathcal{T}\times\mathcal{X}\to\mathbb{R}$ by maximising a generalised empirical likelihood, subject to the restrictions of the sample version of \eqref{UCM}; see \ref{Remark:AiEstimator} in the supplementary file for more details. However, those restrictions are not computable in our context since $T$ is not observable, and the nonparametric estimation of the moment $\mathbb{E}[v(T)]$ for a general function $v(\cdot)$ from contaminated data $\{S_i\}_{i=1}^N$ is challenging and its theoretical properties are difficult to derive, if not impossible. For example, \cite{hall2008estimation} studied the nonparametric estimation of the absolute moment $\mathbb{E}[|T|^{q}]$ with $T$ subject to ordinary smooth error (see the definition \eqref{OrdinaryS}) and found that the theoretical behaviour of the estimator differs depending on $q$: if $q$ is an even integer, the $\sqrt{N}$-consistency is only achievable under a strong condition $\mathbb{E}[T^{2q}]+\mathbb{E}[U^{2q}]<\infty$; if $q$ is an odd integer, the $\sqrt{N}$-consistency is achievable if and only if the distribution of the measurement error is sufficiently ``rough" in terms of the convergence rate of $\phi_U$ to zero in its tails; for $q>0$ not a positive integer, $\sqrt{N}$-consistency is generally impossible. For other forms of $v(T)$ or the involvement of supersmooth error (see the definition \eqref{SuperS}), the consistent nonparametric estimation of $\mathbb{E}[v(T)]$ as well as the corresponding theoretical behaviour are still open problems to our best knowledge.

Thus, to stabilise the estimation of $\pi_0$, we derive another expanding set of equations that can identify $\pi_0$ from the error-contaminated data and avoids estimating $\mathbb{E}[v(T)]$. Specifically, instead of estimating the function $\pi_0(\cdot,\cdot):\mathcal{T}\times\mathcal{X}\to\mathbb{R}$, we turn to estimate its projection $\pi_0(t,\cdot):\mathcal{X}\to\mathbb{R}$ for every \emph{fixed} $t\in\mathcal{T}$, and find that
\begin{equation}
\mathbb{E}\left\{\pi_{0}(t,\bs{X})u(\bs{X})|T=t\right\} =\int_{\mathcal{X}}\frac{f_{T}(t)}{f_{T|X}(t|\bs{x})}u(\bs{x})f_{X|T}(\bs{x}|t)\,d\bs{x}= \mathbb{E}\{u(\bs{X})\}\label{ConditionalMoment1}
\end{equation} holds for any integrable function $u(\bs{X})$. Although the equation~\eqref{ConditionalMoment1} still depends on the unobservable $T$, $\mathbb{E}\left\{\pi_{0}(t,\bs{X})u(\bs{X})|T=t\right\}$ can be estimated using the deconvolution kernel introduced in \eqref{LUCondExp} from the observable $(S,\boldsymbol{X})$. In the following theorem, we show that the corresponding moment condition can identify the function $\pi_{0}(t,\cdot)$ from $(S,\boldsymbol{X})$ for every fixed $t\in\mathcal{T}$.

\begin{thm}\label{Identification}
	Let $L_U(\cdot)$ be the deconvolution kernel function defined in \eqref{KUdef}. For every fixed $t\in\mathcal{T}$ and any integrable function $u(\boldsymbol{X})$, 
	%\begin{equation*}
	%	\lim_{h_0\to 0}\frac{\mathbb{E}\left[\pi(t,\bs{X})u(\bs{X})L_U\{(t-S)/h_0\}\right]}{\mathbb{E}\left[L_U\{(t-S)/h_0\}\right]}=\mathbb{E}\left\{\pi(t,\bs{X})u(\bs{X})|T=t\right\}
	%\end{equation*}	
	%and
	\begin{equation}
		\lim_{h_0\to 0}\frac{\mathbb{E}\left[\pi(t,\bs{X})u(\bs{X})L_U\{(t-S)/h_0\}\right]}{\mathbb{E}\left[L_U\{(t-S)/h_0\}\right]}=\mathbb{E}[u(\bs{X})]\label{ConditionalMoment}
	\end{equation} 
	holds if and only if $\pi(t,\bs{X}) = \pi_{0}(t,\bs{X})$ a.s. 
\end{thm}
The proof is provided in Appendix~A.2 in the supplementary file. Theorem \ref{Identification} suggests a way of estimating the weighting function (i.e. solving a sample analogue of \eqref{ConditionalMoment} for any integrable function $u(\boldsymbol{x})$, where $h_0$ goes to 0 as the sample size tends to infinity). However, this implies solving an infinite number of equations, which is impossible using a finite sample of observations in practice. To overcome this difficulty, we approximate
the infinite-dimensional function space of $u(\boldsymbol{x})$ using a sequence of finite-dimensional
sieves. Specifically, let $u_{K}(\bs{x}):=\big(u_{K,1}(\bs{x}),\ldots,u_{K,K}(\bs{x})\big)^\top$
denote the known basis functions with dimension $K$ (e.g. the power series, B-splines, or trigonometric polynomials). 
The function $u_{K}(\bs{x})$ provides \emph{approximation sieves} that can approximate any suitable functions $u(\bs{x})$ 
 arbitrarily well as $K\to\infty$ (see \citealp{Chen2007} for a discussion on the sieve approximation). Since the 
sieve approximates a subspace of the original function space,  $
\pi_{0}(t,\boldsymbol{X})$ also satisfies  
\begin{equation}
	\lim_{h_0\to 0}\frac{\mathbb{E}\left[\pi(t,\bs{X})u_K(\bs{X})L_U\{(t-S)/h_0\}\right]}{\mathbb{E}\left[L_U\{(t-S)/h_0\}\right]} =\mathbb{E}\{u_K(\bs{X})\}.   \label{sievemoment}
\end{equation}
Equation (\ref{sievemoment}) asymptotically identifies $\pi_{0}(t,\boldsymbol{X})$ as $K\to \infty$. We observe that for any increasing and globally concave function $\rho(v)$, 
\begin{equation}
\pi^\ast(t,\bs{X}):= \rho'\left\{{\lambda^\ast_{t}}^\top u_{K}(\bs{X})\right\}\label{pistar}
\end{equation}
solves (\ref{sievemoment}), where $\rho'(\cdot)$ is the derivative of $\rho(\cdot)$, $\lambda_{t}^\ast :=\argmax_{\lambda\in\mathbb{R}^{K}}G_t^*({\lambda})$ and $G_t^*({\lambda})$ is a strictly concave function defined by
\begin{align*}
G_t^*({\lambda}) := \lim_{h_0\to 0}\frac{\mathbb{E}\left[\rho\{\lambda^\top u_K(\bs{X})\}L_U\{(t-S)/h_0\}\right]}{\mathbb{E}\left[L_U\{(t-S)/h_0\}\right]}-\lambda^\top\mathbb{E}\{u_{K}(\bs{X})\}\,.
\end{align*}
Indeed, by the first-order condition $\nabla G_t^*({\lambda}_t^*)=0$, we see that \eqref{sievemoment} holds with $\pi (t,\bs{X})= \pi^\ast(t,\bs{X})$. The estimator of $\pi_0(t,\bs{X})$ is then expected to be defined as the empirical counterpart of \eqref{pistar}. Therefore, for every fixed $t\in\mathcal{T}$, we propose estimating $\pi_{0}(t,\bs{X})$ by
\begin{equation}
\widehat{\pi}(t,\bs{X}) = \rho'\left\{\widehat{\lambda}_{t}^\top u_{K}(\bs{X})\right\}\label{pihat}
\end{equation}
with $\widehat{\lambda}_{t} = \argmax_{\lambda\in\mathbb{R}^{K}}\widehat{G}_{t}(\lambda)$ and
\begin{align}
\widehat{G}_t(\lambda):=\frac{\sum^{N}_{i=1}\rho\{\lambda^\top u_{K}(\bs{X}_i)\}L_{U}\{(t-S_i)/h_0\}}{\sum^{N}_{i=1}L_{U}\{(t-S_i)/h_0\}}-\lambda^\top\bigg\{\frac{1}{N}\sum^N_{i=1}u_{K}(\bs{X}_i)\bigg\}.\label{hatG}
\end{align} 
Some of the deconvolution kernel $L_{U}\{(t-S_i)/h_0\}$'s may take negative values, making the objective function $\widehat{G}_t(\cdot)$ not strictly concave in a finite sample. However, as $N\to\infty$ and $h_0\to 0$,  $\widehat{G}_t(\cdot)\xrightarrow{p}G^*_t(\cdot)$ and $G^*_t(\cdot)$ is a strictly concave function. Therefore, with probability approaching one, $\widehat{G}_t(\cdot)$ is strictly concave and $\widehat{\lambda}_{t}$ uniquely exists. Remark \ref{rem:TruncateLU}\label{ReplyRemark3} in Section~\ref{sec:ParameterSelection} introduces a way of solving this maximisation problem fast and stably from finite samples.

Our estimator $\widehat{\pi}(t,\boldsymbol{X}_i)$ has a local generalised empirical likelihood interpretation. To see this, Appendix~A.3 shows that $\widehat{\pi}(t,\boldsymbol{X}_i)$ is the dual solution to the following local generalised empirical likelihood maximisation problem: for every fixed $t\in\mathcal{T}$,
\begin{equation}
	\left\{
	\begin{array}
	[c]{cc}
	& \max_{\{\pi_i\}_{i=1}^N} -\frac{\sum_{i=1}^{N}D(\pi_i)L_U(\{t-S_i\}/h_0)}{\sum^N_{i=1}L_U(\{t-S_i\}/h_0)}  \\[2mm]
	& \text{subject to}\ \frac{\sum_{i=1}^{N}\pi_{i}u_{K}(\boldsymbol{X}_{i})L_U(\{t-S_i\}/h_0)}{\sum^N_{i=1}L_U(\{t-S_i\}/h_0)}=  \frac{1}{N}\sum_{i=1}^{N}u_{K}(\boldsymbol{X}%
	_{i})\,,
	\end{array}
	\right.  \label{E:cm1}
\end{equation}
where $D(v)$ is a distance measure from $v$ to 1 for $v\in\mathbb{R}$, which is continuously differentiable and satisfies that $D(1)=0$ and 
$$
\rho(-v) = D\{(D^{'})^{-1}(v)\} - v\cdot (D^{'})^{-1}(v).
$$
Equation \eqref{E:cm1} aims to minimise some distance measure between the desired weight $N^{-1}\pi_i$ and the empirical frequencies $N^{-1}$ locally around a small neighbourhood of $T_i=t$, subject to the sample analogue of the moment restriction~\eqref{sievemoment}.

Since the dual formulation  \eqref{pihat} is equivalent to the primal problem \eqref{E:cm1} and will simplify the following discussions, we shall express the estimator in terms of  $\rho(v)$ in the rest of the discussions. In particular, $\rho(v) = -\exp(-v-1)$ corresponds to exponential tilting \citep{kitamura1997information,imbens1998information,Ai2020}, $\rho(v)=\log(1+v)$ corresponds to the empirical likelihood \citep{owen2001empirical}, $\rho(v)=-(1-v)^2/2$ corresponds to the continuous updating of the generalised method of moments \citep{hansen1982large}, and
$\rho(v) =v -\exp(-v)$ corresponds to the inverse logistic.

Now, replacing $\pi_0(t,\bs{X}_i)$ in \eqref{OracleNonNormalise} with $\widehat{\pi}(t,\bs{X}_i)$, we obtain an estimator of $\mu(t)$:
\begin{equation}
	\widehat{\mu}(t) := \frac{\sum^N_{i=1}\widehat{\pi}(t,\bs{X}_i)Y_iL_U\{(t-S_i)/h\}}{\sum^N_{i=1}L_U\{(t-S_i)/h\}}.\label{CalibrationF}
\end{equation}

\begin{rem}\label{KenndyDR}
	When $T$ is observed without error, \cite{kennedy2017non} propose a doubly-robust estimator for $\mu(t)$ by regressing a pseudo-outcome $$\xi(T,\boldsymbol{X},Y;f_{T|X},m) = \frac{\int f_{T|X}(T|\boldsymbol{x}) f_X(\boldsymbol{x})d\boldsymbol{x}}{f_{T|X}(T|\boldsymbol{X})} \{Y - m(T,\boldsymbol{X})\}+\int m(T,\boldsymbol{x})f_X(\boldsymbol{x})d\boldsymbol{x}$$ onto $T=t$, i.e. $\widehat{\mu}_{DR}(t)=\boldsymbol{g}^\top_{ht}(t) \widehat{\boldsymbol{\beta}}_h(t)$, where $m(T,\boldsymbol{X}):=\mathbb{E}[Y|T,\boldsymbol{X}]$ is the outcome regression function, $\boldsymbol{g}_{ht}(a):=(1,(a-t)/h)^\top$ and
	\begin{align*}
		&\widehat{\boldsymbol{\beta}}_h(t):=\arg\min_{\boldsymbol{\beta}\in\mathbb{R}^2}\frac{1}{Nh}\sum_{i=1}^N L\left(\frac{T_i-t}{h}\right)\left\{\widehat{\xi}(T_i,\boldsymbol{X}_i,Y_i;\widehat{f}_{T|X},\widehat{m})-\boldsymbol{g}_{ht}(T_i)^{\top}\boldsymbol{\beta}\right\}^2,\\
		&\widehat{\xi}(T,\boldsymbol{X},Y;\widehat{f}_{T|X},\widehat{m}) := \frac{N^{-1}\sum_{i=1}^N\widehat{f}_{T|X}(T|\boldsymbol{X}_i)}{\widehat{f}_{T|X}(T|\boldsymbol{X})} \{Y - \widehat{m}(T,\boldsymbol{X})\}+ \frac{1}{N}\sum_{i=1}^N\widehat{m}(T,\boldsymbol{X}_i)
	\end{align*}
	with $L(\cdot)$ being a prespecified kernel function,  $\widehat{f}_{T|X}(\cdot)$ and $\widehat{m}(\cdot)$ are some consistent estimators for $f_{T|X}(\cdot)$ and $m(\cdot)$. \cite{kennedy2017non} showed that $\widehat{\mu}_{DR}(t)$ enjoys double robustness:  (i) when both $f_{T|X}(\cdot)$ and $m(\cdot)$ are consistently estimated and the product of the estimators' local rates of convergence is sufficiently small, $\widehat{\mu}_{DR}(t)$ asymptotically behaves the same as the standard local linear estimator of $\mu(t)$ in \eqref{PC} with known $\pi_0$; (ii) when either $f_{T|X}(\cdot)$ or $m(\cdot)$ is consistently estimated and the other is misspecified, $\widehat{\mu}_{DR}(t)$ is still consistent.

	This idea can be adapted to the our setup with measurement error \eqref{ME} by replacing the standard kernel with the deconvolution one. For example,  we can define a doubly-robust estimator of $\mu(t)$ by
	\begin{align*}
		\widehat{\mu}_{DK}(t):=\frac{\sum^N_{i=1}\widehat{\xi}_{DK}(t,\boldsymbol{X}_i,Y_i;\widehat{\pi},\widehat{m}_{DK}) L_U\{(t-S_i)/h\}}{\sum^N_{i=1}L_U\{(t-S_i)/h\}},
	\end{align*} 
	where
	$$\widehat{\xi}_{DK}(t,\boldsymbol{X},Y;\widehat{\pi},\widehat{m}_{DK}) = \widehat{\pi}(t,\boldsymbol{X}) \{Y - \widehat{m}_{DK}(t,\boldsymbol{X})\}+ \frac{1}{N}\sum_{i=1}^N\widehat{m}_{DK}(t,\boldsymbol{X}_i)$$
	and $\widehat{m}_{DK}(\cdot)$ is some consistent estimator of $m(\cdot)$. Comparing to our proposed estimator  $\widehat{\mu}(t)$ in \eqref{CalibrationF},  $\widehat{\mu}_{DK}(t)$ requires additionally a consistent estimator of $m(\cdot)$ from the error-contaminated data $(S_i,\bs{X}_i,Y_i)_{i=1}^N$. Establishing practical estimators (with tuning parameter techniques) and the corresponding theoretical results for this method is beyond the scope of this paper and will be resolved in future work.
\end{rem}

%\whcomment{What do we want to do with this then? From Haoze's work, I feel it may not be hard to derive the theory for this estimator based on what we have had. But of course, we do not want to make this the main estimator of this paper, since it will change a lot the whole paper and takes a lot of time to redo the numerical works. So, how do we put this? Maybe a shorter remark here saying that, once $\pi_0$ is estiamted robustly by our proposed localized GEL estimator $\widehat{\pi}$, a nonparametric doubly-robust estimator of $\mu$ under Kennedy's framework is also possible. Then discuss that it requires an addition deconvolution kernel regression estimator of the outcome regression function $m(t,\bs{x})$, but it is beneficial in a weaker condition on the tuning parameters. Then we put the details and maybe the theoretical results of this estimator in the appendix?}

\section{Large Sample Properties}\label{sec:asymptotic}
In this section, we establish the $L_\infty$ and $L_2$ convergence rates of $\widehat{\pi}(t,\cdot)$ for every fixed $t\in\mathcal{T}$. We then investigate the asymptotic behaviour of the proposed ADRF estimator $\widehat{\mu}(t)$. Note that $\widehat{\mu}(t)$ (resp. $\widehat{\pi}(t,\cdot)$) is a nonparametric estimator and that its asymptotic behaviour is affected by the asymptotic bias and variance, which are respectively defined as the expectation and variance of the limiting distribution of $\widehat{\mu}(t)-\mu(t)$ (resp. $\widehat{\pi}(t,\cdot)-\pi_0(t,\cdot)$). Based on \eqref{LUCondExp}, we will show that the asymptotic biases of the two estimators are the same as their counterparts in the error-free case.  That is, they depend on the smoothness of $\pi_0$, $\mu$, and the density of $T$, and the approximation error based on the sieve basis $u_K$. In particular, the following conditions are required:

\begin{assumption}\label{as:KernelOrder}
	The kernel function $L(\cdot)$ is an even function such that $\int_{-\infty}^\infty L(u)\,du =1 $ and has finite moments of order 3.
\end{assumption}

\begin{assumption} We assume \label{as:positivity}
	\begin{enumerate}
	 \item [(i)]   the support $\mathcal{X}$ of $\boldsymbol{X}$ is a
	compact subset of $\mathbb{R}^{r}$. The support $\mathcal{T}$ of the
	treatment variable $T$ is a compact subset of $\mathbb{R}$. 
	\item [(ii)] (Strict Positivity)  there exist
	a positive constant $\eta_{min}$ such that 
	$f_{T|X}(t|\boldsymbol{x})\geq \eta_{\min}>0$, for all  $\boldsymbol{x} \in\mathcal{X}$.
	\end{enumerate}
\end{assumption}

\begin{assumption}\label{as:smooth_densities}
 (i) The densities $f_T(t)$, $f_{T|X}(t|\boldsymbol{X})$ and $f_{T|Y,X}(t|Y,\boldsymbol{X})$ are third-order continuously differentiable w.r.t. $t$ almost surely. (ii)  The derivatives of  $f_{T|X}(t|\boldsymbol{X})$ and  $f_{T|Y,X}(t|Y,\boldsymbol{X})$, denoted by $\{\partial_t^df_{T|X}(t|\boldsymbol{X}),\ \partial_t^df_{T|Y,X}(t|Y,\boldsymbol{X})\ \text{for} \ d=0,1,2,3\}$, are integrable almost surely in $t$.
\end{assumption}

\begin{assumption}
	\label{as:smooth_pi} For every $t\in\mathcal{T}$, (i) the function $\pi_0(t,\bs{x})$ is $s$-times continuously differentiable w.r.t. $\boldsymbol{x}\in\mathcal{X}$, where $s>r/2$ is an integer; (ii)  there exist $\lambda_{t}\in\mathbb{R}
	^{K}$ and a positive constant $\alpha>0$ such that
	$\sup_{\boldsymbol{x}\in \mathcal{X}}\left\vert
	(\rho^{\prime})^{-1}\left\{\pi_{0}(t,\boldsymbol{x})\right\}  -\lambda_{t}^\top u_{K}(\boldsymbol{x})\right\vert
	=O(K^{-\alpha})$.
\end{assumption}

\begin{assumption}
	\label{as:u&v}(i) For every $K$, the eigenvalues of $\mathbb{E}
	\left[   u_{K}(\boldsymbol{X})u_{K}(\boldsymbol{X})^{\top}|T=t\right] $
	are bounded away from zero and infinity, and twice differentiable w.r.t. $t$ for $t\in\mathcal{T}$. (ii)  There is
	a  sequence of constants $\zeta(K)$ satisfying $\sup_{ \boldsymbol{x}\in\mathcal{X}}\Vert u_{K}(\boldsymbol{x}%
	)\Vert\leq\zeta(K)$, such that $\zeta(K)\{K^{-\alpha}+h_0^2+h^2\}\rightarrow0$
	as $N\rightarrow\infty$, where $\|\cdot\|$ denotes the  Euclidean norm.
\end{assumption}

\begin{assumption}
	\label{as:smooth_m} For every $t\in\mathcal{T}$, there exist $\gamma_{t}\in\mathbb{R}
	^{K}$ and a positive constant $\ell>0$ such that
	$\sup_{\boldsymbol{x}\in \mathcal{X}}\left\vert
	m(t,\boldsymbol{x}) -\gamma_{t}^\top u_{K}(\boldsymbol{x})\right\vert=O(K^{-\ell})$, where $m(t,\boldsymbol{x})=\mathbb{E}[Y|T=t,\boldsymbol{X}=\boldsymbol{x}]$.
\end{assumption}

\begin{assumption}\label{as:Y_bounds}
	$R_1^{2+\delta}(t):=\mathbb{E}\big[|\pi_0(t,\boldsymbol{X})Y - \mu(t)|^{2+\delta}|T=t\big]$, $R_2^{2+\delta}(t):=\mathbb{E}\big[|\pi_0(t,\boldsymbol{X})m(t,\bs{X})- \mu(t)\big|^{2+\delta}|T=t\big]$ and $R_3^{2+\delta}(t):=\mathbb{E}\big[|\pi_0(t,\boldsymbol{X})\{Y-m(t,\bs{X})\}\big|^{2+\delta}|T=t\big]$ are bounded for some $\delta>0$, for all $t\in\mathcal{T}$.%, where $m(t,\bs{x}):= \mathbb{E}(Y|T=t,\bs{X}=\bs{x})$.
\end{assumption}

Assumption \ref{as:positivity} (i) restricts the covariates and the treatment to be bounded. This condition is commonly imposed in the nonparametric regression literature. Assumption~\ref{as:positivity} (i) can be relaxed if we restrict the tail distributions of $\boldsymbol{X}$ and $T$. For example, \citet[Assumption 3]{chen2008semiparametric}  allowed the support of $\boldsymbol{X}$ to be the entire Euclidean space but imposed $\int_{\mathbb{R}^r} (1+|\boldsymbol{x}|^2)^{\omega}f_X(\boldsymbol{x})d\boldsymbol{x}<\infty$  for some $\omega>0$.
	
Assumption \ref{as:positivity} (ii) is a strict positivity condition requires every subject having certain chance of receiving every treatment level regardless of covariates. This condition is also imposed in a large body of literature  in the absence of measurement error (see e.g.~\citealp[Assumption 2]{kennedy2017non} and \citealp[Assumption 3]{d2021overlap}), particularly when no restrictions are imposed on the potential outcome distribution.
This condition can be relaxed if other smoothness  conditions are imposed on the potential outcome distribution \citep{ma2020robust}, or if different target parameters are considered; for example, \cite{munoz2012population} studied the  estimation of a \emph{stochastic} intervention causal parameter, defined by  $\mathbb{E}[\mathbb{E}[Y|T+a(\boldsymbol{X}),\boldsymbol{X}]]$, based on a weaker positivity condition, i.e. $\sup_{t\in\mathcal{T}}\{f_{T|X}(t-a(\boldsymbol{X})|\boldsymbol{X})/f_{T|X}(t|\boldsymbol{X})\}<\infty$ a.e., where $a(\boldsymbol{X})$ is a user specified  intervention function. \cite{diaz2013targeted}  studied the estimation of a conditional causal dose-response curve defined by $\mathbb{E}[\mathbb{E}[Y|T=t,\boldsymbol{X}]|\boldsymbol{Z}]$, where $\boldsymbol{Z}\subset \boldsymbol{X}$ is a subset of observed covariates, based on a  weaker positivity condition, i.e. $\sup_{t\in\mathcal{T}}\{b(t,\boldsymbol{Z})/f_{T|X}(t|\boldsymbol{X})\}<\infty$ a.e., for a user specified weight function $b(t,\boldsymbol{Z})$.   Although Assumption \ref{as:positivity} (ii) is not the  mildest condition in the literature, we maintain it throughout this paper owing to its technical benefits, especially in the presence of measurement error.\label{ReplyPositivity}
	
Assumption~\ref{as:smooth_densities}\label{ReplySmoothness} includes smoothness conditions required for nonparametric estimation.     Under Assumption \ref{UnconfoundAssump}, the parameter of interest $\mu(t)=\mathbb{E}[Y^*(t)]$ can be also written as $\mu(t) =\mathbb{E}\left[\pi_0(t,\boldsymbol{X})Y|T=t\right]
%=  \int \pi_0(t,\boldsymbol{x}) \frac{f_{T|Y,X}(t|y,\boldsymbol{x})}{f_{T}(t)} \cdot y f_{X,Y}(\boldsymbol{x},y)d\boldsymbol{x}dy\\
=\mathbb{E}\left[f_{T|Y,X}(t|Y,\boldsymbol{X})Y/f_{T|X}(t|\boldsymbol{X}) \right]$.
Note that Assumption \ref{as:smooth_densities} (i) implies that $t \mapsto  f_{T|Y,X}(t|Y,\boldsymbol{X})/f_{T|X}(t|\boldsymbol{X})$ is third-order continuously differentiable almost surely. Furthermore, using Leibniz integral rule and Assumption \ref{as:smooth_densities} (ii), we have that the target parameter $t\to \mu(t)$ is third-order continuously differentiable.
	
Assumption \ref{as:smooth_pi} (i)\label{ReplyMildness}  is used to control the complexity (measured by the uniform entropy integral) of the function class $\{\pi_0(t,\boldsymbol{x}),\boldsymbol{x}\in\mathcal{X}\}$ such that it forms a Donsker class and the empirical process theory can be applied \citep[Corollary 2.7.2]{van1996weak}. Despite of its stringency,  the smoothness condition of this type is commonly adopted in the literature of nonparametric inference, see  \citet[Assumption 4 (i)]{chen2008semiparametric} and \citet[Condition E.1.7]{fan2021optimal}.

Assumption~\ref{as:smooth_pi} (ii) requires the sieve approximation error of $\rho^{\prime-1}\left\{\pi_{0}(t,\boldsymbol{x})\right\}$ to shrink at a polynomial rate. This condition is satisfied for a variety of sieve basis functions. For example, it can be satisfied with $\alpha=+\infty$ if $\bs{X}$ is discrete, and with $\alpha = s/r$ if $\bs{X}$ is continuous and $u_K(\bs{x})$ is a power series or a B-spline, where $s$ is the smoothness of the approximand and $r$ is the dimension of $\bs{X}$. Assumption~\ref{as:smooth_m} imposes a similar sieve approximation error for $m(t,\bs{x})$.

Assumption \ref{as:u&v} (i) rules out near multicollinearity in the approximating basis functions, which is common in the sieve regression literature. Assumption~\ref{as:u&v}~(ii) is satisfied with $\zeta(K)=O(K)$ if $u_K(\bs{x})$ is a power series and with $\zeta(K)=O(\sqrt{K})$ if $u_K$ is a B-spline \citep{Newey97}. Assumption \ref{as:Y_bounds} imposes the boundedness conditions on the moment of the response variable, which are also standard in the errors-in-variables problem (e.g.~\citealp{Fan1993, delaigle2009design}). This condition is needed for deriving the asymptotic distribution of the proposed estimator by applying the Lyapunov central limit theorem.

Depending on the type of the distribution of $U$ and decaying rates of $h_0$ and $h$, the asymptotic variance of our estimator differs. This is different from the error-free case. We consider two types of $U$: the ordinary smooth case and supersmooth case, which are standard in the literature of errors-in-variables problem (see e.g. \citealp{Fan1993}, \citealp{delaigle2009design}, and \citealp{Meister2009}, among others, for more details).
		
An ordinary smooth error of order $\beta \geq 1$ satisfies
\begin{equation}
\lim_{t\rightarrow\infty} t^\beta \phi_U(t) = c \quad \text{and} \quad \lim_{t\rightarrow\infty} t^{\beta+1}\phi_U^{(1)}(t) = -c\beta\,,\label{OrdinaryS}
\end{equation}
for some constant $c>0$. A supersmooth error of order $\beta \geq 1$ satisfies
\begin{equation}\label{SuperS}
d_0|t|^{\beta_0}\exp(-|t|^\beta/\gamma) \leq |\phi_U(t)|\leq d_1|t|^{\beta_1}\exp(-|t|^\beta/\gamma) \quad \text{as} \quad |t| \rightarrow \infty\,,
\end{equation}
for some positive constants $d_0,d_1,\gamma$ and some constants $\beta_0$ and $\beta_1$.\label{ReplyOrderU}
Examples of ordinary smooth errors include Laplace errors, Gamma errors, and their convolutions. Cauchy errors, Gaussian errors, and their convolutions are supersmooth errors. The order $\beta$ describes the decaying rate of the characteristic function $\phi_U(t)$ as $t\rightarrow\infty$, which corresponds to the smoothness of the error distribution (e.g.~$\beta=1$ for Cauchy distribution, $\beta=2$ for Laplace and Gaussian distribution and for Gamma distribution, it relates to both the shape and the scale parameters).

Since in the inverse Fourier transform representation \eqref{eq:fT_InverseF}, division by $\phi_U$ appears,  it is natural to expect  better estimation results for a larger $|\phi_U|$ (i.e. a smaller $\beta$); indeed, it is found in the literature (see e.g. \citealp{fan1991optimal}, \citealp{Fan1993}, \citealp{delaigle2009design}, and \citealp{Meister2009}, among others) that for both the ordinary smooth and supersmooth cases, the higher the order $\beta$ is, the harder the deconvolution will be, i.e. the slower the variance of a deconvolution kernel estimator converges. This is an intrinsic difficulty to the nonparametric estimation with errors in variables  (\citealp{Fan1993}, \citealp{Carroll2006}).

Such an influence will be seen in the convergence rate of our estimator in the following theorems.
%Lemmas~\ref{OrdinaryVar} and \ref{SuperVar} in Appendix B, adapted from \cite{Fan1993} and \cite{delaigle2009design}, show the extent to which the distribution of $U$ influences the asymptotic variance of a deconvolution kernel estimator using $L_U$. 
Depending on the type of the distribution of $U$, we need the following different conditions on $L$ to derive the asymptotic variance:\\[.2cm]
\noindent
{{\bf Assumption O} (Ordinary Smooth Case):} $\|\phi_L\|_{\infty}<\infty$, $\int_{-\infty}^\infty |t|^{\beta+1} \{|\phi_L(t)|+|\partial_t\phi_L(t)|\}\,dt<\infty$ and $\int_{-\infty}^\infty |t^{\beta}\phi_L(t)|^2\,dt <\infty $.
\\[.2cm]
\noindent
{{\bf Assumption S} (Supersmooth Case):} $\phi_L(t)$ is support on $[-1,1]$ and bounded.

These assumptions concern the prespecified kernel function and can be satisfied easily. For example, the one whose Fourier transform is $\phi_{L}(u) = (1-u^2)^3\cdot \mathbbm{1}_{[-1,1]}(u)$ satisfies these conditions (e.g.~\citealp{Fan1993} and \citealp{delaigle2009design}). In the following two sections, we establish the large sample properties of $\widehat{\pi}(t,\cdot)$ and $\widehat{\mu}(t)$ under the two types of $U$.

\subsection{Asymptotics for the Ordinary Smooth Error}
To establish the large sample properties of $\widehat{\mu}(t)$, we first show that the estimated weight
function $\widehat{\pi}(t,\cdot)$ is consistent and compute its
convergence rates under both the $L_{\infty}$ norm and the $L_{2}$ norm.

\begin{thm}\label{thm:hatp-p}
	Suppose that the error $U$ is ordinary smooth of order $\beta$ satisfying \eqref{OrdinaryS} and that Assumption O holds.	Under Assumptions \ref{as:KernelOrder}--\ref{as:u&v} and $\zeta(K)\sqrt{K\Big/\left(Nh_0^{1+2\beta}\right)}\rightarrow0$ as $N\rightarrow \infty$,
		 for every fixed $t\in\mathcal{T}$, then  
		\begin{align*}
		&\sup_{\bs{x}\in\mathcal{X}}|\widehat{\pi}(t,\bs{x})-\pi_0(t,\bs{x})|=O_p\left(\zeta(K)\left\{K^{-\alpha}+h_0^2\right\}+\zeta(K)\sqrt{\frac{K}{Nh_0^{1+2\beta}}}\right), \\
		&\int_{\mathcal{X}}|\widehat{\pi}(t,\bs{x})-\pi_0(t,\bs{x})|^2dF_X(\bs{x})=O_p\left(\{K^{-2\alpha}+h_0^4\}+{\frac{K}{Nh_0^{1+2\beta}}}\right), \\
		&\frac{1}{N}\sum_{i=1}^N|\widehat{\pi}(t,\bs{X}_i)-\pi_0(t,\bs{X}_i)|^2=O_p\left(\{K^{-2\alpha}+h_0^4\}+{\frac{K}{Nh_0^{1+2\beta}}}\right). 
		\end{align*} 
	\end{thm}	
The proof of Theorem \ref{thm:hatp-p} is presented in Appendix~C. 
The first part of the rates, $\zeta(K)\left\{K^{-\alpha}+h_0^2\right\}$ and $\{K^{-2\alpha}+h_0^4\}$, are the rates of the asymptotic bias. $\zeta(K)\sqrt{K/Nh_0^{1+2\beta}}$ and $K/Nh_0^{1+2\beta}$ correspond to the asymptotic variance.

We next establish the asymptotic linear expansion and asymptotic normality of $\widehat{\mu}(t) - \mu(t)$. To aid the presentation, we define the following quantities. For $i=1,\ldots,N$, $\eta_{h,h_0}(S_i,\boldsymbol{X}_i,Y_i;t):=\phi_{h}(S_i,\boldsymbol{X}_i,Y_i;t)+\psi_{h_0}(S_i,\boldsymbol{X}_i,Y_i;t)$, where
\begin{align*}
	\phi_{h}&(S_i,\boldsymbol{X}_i,Y_i;t):=\big[\pi_0(t,\boldsymbol{X}_i)Y_iL_{U,h}(t-S_i)-\mathbb{E}\{\pi_0(t,\boldsymbol{X})YL_{U,h}(t-S)\}\big]\\
    &-\mu(t)\big[L_{U,h}(t-S_i)-\mathbb{E}\{L_{U,h}(t-S)\}\big],\\
	\psi_{h_0}&(S_i,\boldsymbol{X}_i,Y_i;t):=\mu(t)\big[L_{U,h_0}(t-S_i)-\mathbb{E}\{L_{U,h_0}(t-S)\}\big]\\
	&- \big[m(t,\bs{X}_i)\pi_0(t,\bs{X}_i) L_{U,h_0}\left(t-S_i\right)-\mathbb{E}\{m(t,\bs{X})\pi_0(t,\bs{X})L_{U,h_0}\left(t-S\right)\}\big],
\end{align*}
with $L_{U,h}(v):=h^{-1}L_U(v/h)$. The population mean of both $\phi_{h}$ and $\psi_{h_0}$ are zero. Let $``\ast"$ denote the convolution operator, we define
\begin{align*}
	V_{j}:=&f_T^{-2}(t) (R_j^2f_T)\ast f_U(t)\cdot C\,,\ \text{for} \ j=1,2\,,
\end{align*}
where $C:=\int^\infty_{-\infty} J^2(v)\,dv = (2\pi c^2)^{-1}\int |w|^{2\beta}\phi_L^2(w)\,dw$, with $c$ defined in \eqref{OrdinaryS}, $J(v):= (2\pi c)^{-1}\int^\infty_{-\infty}\exp(-iwv)\phi_L(w)w^\beta\,dw$ and $R_1^2, R_2^2$ defined in Assumption~\ref{as:Y_bounds}. 
Moreover, let 
$(R_1R_2)(t) := \mathbb{E}\big[\{\pi_0(t,\boldsymbol{X})Y - \mu(t)\}\{\mu(t)-\pi_0(t,\boldsymbol{X})m(t,\boldsymbol{X})\}|T=t\big]$ and  $v_{h}(t):=\mathbb{E}\{L^2_{U,h}(t-S)\}$. 

\begin{thm}\label{thm:normality}
Suppose that the error $U$ is ordinary smooth of order $\beta$ satisfying \eqref{OrdinaryS} and that Assumption O and Assumptions \ref{UnconfoundAssump}--\ref{as:Y_bounds} as well as the following condition hold:
$$\frac{(K^{-\ell}+h_0^2)\cdot (K^{-\alpha}+h_0^2)}{h^2}+\frac{(h\wedge h_0)^{1/2+\beta}}{h_0^{1+2\beta}}\frac{K}{\sqrt{N}}\to 0,$$	 
where $(h\wedge h_0) = h\mathbbm{1}\{h=O(h_0)\}+h_0\mathbbm{1}\{h_0=o(h)\}$. Then, for every fixed $t\in\mathcal{T}$,
\begin{equation}
\begin{aligned}
	\hspace{-.3cm}\widehat{\mu}(t)- \mu(t)=&\label{LinearExpressionOrdinary}\frac{\kappa_{21}}{2}\bigg[\frac{f_T(t)\Phi_1(t)-{\mu(t)\partial_t^2f_T(t)}}{f_T(t)}\bigg]\cdot h^2 +o(h^2)\\
	+&\frac{\kappa_{21}}{2} \bigg[\frac{\mu(t)\partial_{t}^2f_T(t) - f_T(t)\Phi_2(t)}{f_T(t)}\bigg]\cdot h_0^2 + o(h_0^2)\\
+& \sum_{i=1}^N\frac{\eta_{h,h_0}(S_i,\boldsymbol{X}_i,Y_i;t)}{N\cdot f_T(t)} + o_P\bigg\{\frac{1}{\sqrt{N(h\wedge h_0)^{1+2\beta}}}\bigg\}\,,
\end{aligned}
\end{equation}
where $\kappa_{ij}:=\int u^i L^j(u)du$, $\Phi_1(t):= \mathbb{E}\big[\{Y\partial_t^2f_{T|Y,\bs{X}}(t|Y,\bs{X})\}/\{f_{T|\bs{X}}(t|\bs{X})\}\big]$, and $\Phi_2(t):= \mathbb{E}\big[\{m(t,\bs{X})\partial_t^2f_{T|\bs{X}}(t|\bs{X})\}/\{f_{T|\bs{X}}(t|\bs{X})\}\big]$. Furthermore,
\begin{enumerate}[a)]
	\item if $h = o(h_0)$, then $\sqrt{h^{1+2\beta}/N}\sum_{i=1}^N\eta_{h,h_0}(S_i,\boldsymbol{X}_i,Y_i;t)/f_T(t)\overset{d}{\to}N(0,V_1);$
	\item if $h_0 = o(h)$, then $\sqrt{h_0^{1+2\beta}/N}\sum_{i=1}^N \eta_{h,h_0}(S_i,\boldsymbol{X}_i,Y_i;t)/f_T(t)\overset{d}{\to}N(0,V_2);$
	\item if $h_0 = \tilde{c}h$ for a constant $\tilde{c}>0$, then 
	$\sqrt{h^{1+2\beta}/N}\sum_{i=1}^N \eta_{h,h_0}(S_i,\boldsymbol{X}_i,Y_i;t)/f_T(t)\overset{d}{\to}N(0,V_3)\,,$ where 
	\begin{align*}
		V_3:=& \frac{(R_1^2f_T)\ast f_U(t)}{f^2_T(t)} \cdot \int_{-\infty}^\infty J^2(v)\,dv + \frac{(R_2^2f_T)\ast f_U(t)}{\tilde{c}^{(2+2\beta)}f^2_T(t)}\cdot \int_{-\infty}^\infty J^2(v/\tilde{c})\,dv\\
		&+ \frac{2\{(R_1R_2)f_T\}\ast f_U(t)}{\tilde{c}^{(1+\beta)}f_T^2(t)}\cdot \int_{-\infty}^\infty J(v)J(v/\tilde{c})\,dv\,.
	\end{align*}
In particular, when $\tilde{c}=1$, $V_3$ reduces to $f_T^{-2}(t) (R_3^2f_T)\ast f_U(t)\cdot C$ with $R_3^2$ defined in Assumption~\ref{as:Y_bounds}.
	
\end{enumerate}
\end{thm}
%\whcomment{Do we want to leave the explicit expressions $V_1,V_2$ and $V_3$ in the main text or in the appendix? I prefere to put them in the appendix. Here we can write similarly to the supersmooth case:
%Furthermore,
%$$
%[var\{\eta_{h,h_0}(S,\bs{X},Y;t)\}]^{-1/2}\frac{1}{\sqrt{N}}\sum^N_{i=1}\eta_{h,h_0}(S_i,\bs{X}_i,Y_i;t) \overset{D}{\to} N(0,1)\,,
%$$
%where 
%$$
%var\{\eta_{h,h_0}(S,\bs{X},Y;t)\} \asymp (h \wedge h_0)^{-(1+2\beta)}%\,.
%$$
%}
The proof of Theorem \ref{thm:normality} is presented in Appendix~D. %In our case, the estimation involves an $r$-dimensional variable $\bs{X}$, where $r\geq 1$. The curse of dimensionality presents in $K^{-\alpha}=K^{-s/r}$ in the bias arising from $\widehat{\pi}(t,\cdot)$, where $s$ is the smoothness of $\pi_0(t,\boldsymbol{x})$ w.r.t.~$\boldsymbol{x}$.
From the theorem, we see that as long as $\pi_0(t,\cdot)$  and $m(t,\cdot)$ are sufficiently smooth or $K$ grows sufficiently fast, and $h_0$ decays fast enough, so that $(K^{-\ell}+h_0^2)\cdot (K^{-\alpha}+h_0^2) = o(h^2)$, the error arising from the sieve approximation is asymptotically negligible. For example, using the usual trade-off between the squared bias and variance, $\widehat{\mu}(t)-\mu(t)$ achieves the optimal convergence rate, $N^{-2/(2\beta+5)}$, if $h_0\asymp h \asymp N^{-1/(2\beta+5)}$. In such a case, we require $K=o(h^{-2})$, $\alpha+\ell>1$ and $\alpha>1/2$ if spline basis is used and $\alpha>1$ if a power series is used (The detailed derivation can be found in Appendix~A.5).

The convergence rate $N^{-2/(2\beta+5)}$ above is optimal for all possible nonparametric regression estimators when the regressors are measured with ordinary smooth errors showed in \cite{Fan1993}. Note that for error-free local constant estimator, the convergence rates of the asymptotic bias and variance are $h^2$ and $(Nh)^{-1/2}$, respectively (see e.g.~\citealp{Fan1996}). Our proposed estimator $\widehat{\mu}(t)$ has the same rate of asymptotic bias as that in the error-free case, but the asymptotic variance is degenerated by $h^{-\beta}$, owing to the ordinary smoothness of the error distribution.\label{ReplyOrdinaryRate}

In addition to asymptotic normality, we provide in \eqref{LinearExpressionOrdinary} the asymptotic linear expansion of $\widehat{\mu}(t) - \mu(t)$, which can help conduct statistical inference. It is known in the literature on measurement error (see e.g.~\citealp{Delaigle2015} Appendix C) \label{ReplyVarianceEstimation} that the closed-form asymptotic variances $V_1,V_2$ and $V_3$ are difficult to estimate. However, using our linear expansion in \eqref{LinearExpressionOrdinary}, to estimate the asymptotic variance, we only need consistent estimators of $\phi_{h}(S,\boldsymbol{X},Y;t)$ and $\psi_{h_0}(S,\boldsymbol{X},Y;t)$. For example, $\pi_0(t,\bs{X})$ and $\mu(t)$ can be estimated respectively using our $\widehat{\pi}(t,\bs{X})$ and $\widehat{\mu}(t)$, and $\mathbb{E}[Y|T=t,\boldsymbol{X}]$ can be estimated using \citeauthor{liang2000asymptotic}'s (\citeyear{liang2000asymptotic}) method. Then, we can construct a pointwise confidence interval for $\mu(t)$ using the undersmoothing technique (see Appendix~A.4 for the detailed method and some simulation results). Other confidence intervals based on bias-correction \citep[see e.g.][]{calonico2018effect,takatsu2022debiased} are also possible but require a better estimation of the asymptotic bias and corresponding adjustments of the variance estimation with theoretical justification. That is beyond the scope of this paper and will be resolved in future work.

\subsection{Asymptotics for the Supersmooth Error}
The next two theorems establish the asymptotic properties of our estimator for the supersmooth case.
\begin{thm}\label{thm:hatp-p_SS}
	Suppose that the error $U$ is supersmooth of order $\beta$ satisfying \eqref{SuperS} and Assumption S holds. Under Assumptions \ref{as:KernelOrder}--\ref{as:u&v} and $\zeta^2(K)K\cdot(Nh_0)^{-1}\cdot\exp(2h_0^{-\beta}/\gamma)\rightarrow 0$ as $N\rightarrow \infty$,
	for every fixed $t\in\mathcal{T}$, then  
	\begin{align*}
	&\sup_{\bs{x}\in\mathcal{X}}|\widehat{\pi}(t,\bs{x})-\pi_0(t,\bs{x})|=O_p\left(\zeta(K)\cdot\left[\left\{K^{-\alpha}+h_0^2\right\}+\frac{\exp\left(h_0^{-\beta}/\gamma\right)}{\sqrt{h_0}}\cdot\sqrt{\frac{K}{N}}\right]\right), \\
	&\int_{\mathcal{X}}|\widehat{\pi}(t,\bs{x})-\pi_0(t,\bs{x})|^2dF_X(\bs{x})=O_p\left(\{K^{-2\alpha}+h_0^4\}+\frac{\exp\left(2h_0^{-\beta}/\gamma\right)}{h_0}\cdot {\frac{K}{N}}\right), \\
	&\frac{1}{N}\sum_{i=1}^N|\widehat{\pi}(t,\bs{X}_i)-\pi_0(t,\bs{X}_i)|^2=O_p\left(\{K^{-2\alpha}+h_0^4\}+\frac{\exp\left(2h_0^{-\beta}/\gamma\right)}{h_0}\cdot {\frac{K}{N}}\right). 
	\end{align*} 
\end{thm}	
The proof of Theorem~\ref{thm:hatp-p_SS} is presented in Appendix~C. Comparing these results with those in Theorem~\ref{thm:hatp-p}, the asymptotic bias is the same as that in the ordinary smooth case. The rate of the asymptotic variance, however, becomes much slower, which is expected in the errors-in-variables context; see~\cite{Fan1993} and \cite{delaigle2009design}.

\begin{thm}\label{thm:normality_SS}
Suppose that the error $U$ is supersmooth of order $\beta$ satisfying \eqref{SuperS} and that Assumption S and Assumptions \ref{UnconfoundAssump}--\ref{as:Y_bounds} hold. Letting $e(h) := h^{1/2}\exp(-h^{-\beta}/\gamma)$, we have $v_{h}(t)=\mathbb{E}\{L^2_{U,h}(t-S)\}=O\{e(h)^{-2}\}$. If, as $h\rightarrow 0$, $v_{h}(t)\rightarrow\infty$ and
	$$
	\frac{(K^{-\ell}+h_0^2) \cdot (K^{-\alpha}+h_0^2)}{h^2}
	+\frac{K}{\{e(h) \wedge e(h_0)\}\sqrt{N}}\to 0 \ \text{as} \ N\rightarrow\infty\,,
	$$
	then, for every fixed $t\in\mathcal{T}$,
	\begin{equation}
		\begin{aligned}
			\widehat{\mu}(t)- \mu(t)=&\frac{\kappa_{21}}{2}\bigg[\frac{f_T(t)\Phi_1(t)-{\mu(t)\partial_t^2f_T(t)}}{f_T(t)}\bigg]\cdot h^2 +o(h^2)\\
			+&\frac{\kappa_{21}}{2} \bigg[\frac{\mu(t)\partial_{t}^2f_T(t) - f_T(t)\Phi_2(t)}{f_T(t)}\bigg]\cdot h_0^2 + o(h_0^2)\\
			&+ \sum_{i=1}^N\frac{\eta_{h,h_0}(S_i,\boldsymbol{X}_i,Y_i;t)}{N\cdot f_T(t)}\cdot \{1+o_P(1)\}\,,\label{LinearExpressionSS}
		\end{aligned}
	\end{equation}
	where $\kappa_{21},\Phi_1(t),\Phi_2(t)$ and $\eta_{h,h_0}$ are defined as those in Theorem~\ref{thm:normality} and
	$$
	\{Nf_T(t)\}^{-1} \sum_{i=1}^N \eta_{h,h_0}(S_i,\boldsymbol{X}_i, Y_i; t)=O_p\{N^{-1/2}\{e(h)\wedge e(h_0)\}^{-1}\}\,.
	$$
	Moreover,
	%\begin{enumerate}[a)]
	%	\item 
		if $v_h(t)\geq d_1 f_S(t)h^{d_3}\exp(2h^{-\beta}/\gamma - d_2h^{-d_4\beta})$ for some constants $d_1,d_2>0$, $1>d_4>0$ and $d_3$, we have
		$$
		[\text{var}\{\eta_{h,h_0}(S_i,\boldsymbol{X}_i,Y_i;t)\}]^{-1/2}\cdot\frac{1}{\sqrt{N}}\sum_{i=1}^N\big\{{\eta_{h,h_0}(S_i,\boldsymbol{X}_i,Y_i;t)}\big\} \overset{D}{\to} N(0,1)\,.
		$$
		%\item The optimal convergence rate of $\widehat{\mu}(t)$ is $\left( \log N\right)^{-2/\beta}$, which is obtained when $h \wedge h_0 = h=d(\log N)^{-1/\beta}$ for a constant $d>(2/\gamma)^{1/\beta}$. In this case, the rate of variance, $N^{-1}\{e(h)\wedge e(h_0)\}^{-2}=o(h^4)$, is negligible.
		%\item if $h \wedge h_0 =d(\log N)^{-1/\beta}$, where $d>(2/\gamma)^{1/\beta}$, then, $N^{-1/2}\{e(h)\wedge e(h_0)\}^{-1}=o(h^2)$. Then
		%\begin{enumerate}
		%	\item if $h \wedge h_0 = h$, then, for every fixed $t\in\mathcal{T}$,
		%	\begin{align*}
		%		&\widehat{\mu}(t)- \mu(t)=\frac{\kappa_{21}}{2}\bigg[\frac{f_T(t)\Phi(t)-\mu(t)\partial_t^2f_T(t)}{f_T(t)}\bigg]\cdot {\color{red}d^2\left( \log N\right)^{-2/\beta}}\cdot\{1+o_P(1)\}\,;
		%	\end{align*}
		%	\item if $h \wedge h_0 = h_0$, then, for every fixed $t\in\mathcal{T}$,
		%	\begin{align*}
		%		&\widehat{\mu}(t)- \mu(t)=\frac{\kappa_{21}}{2}\bigg[\frac{f_T(t)\Phi(t)-\mu(t)\partial_t^2f_T(t)}{f_T(t)}\bigg]\cdot h^2\cdot\{1+o_P(1)\}\,.
		%	\end{align*}
		%\end{enumerate}
	%\end{enumerate}
	
\end{thm}
The proof of Theorem~\ref{thm:normality_SS} is presented in Appendix~D. As in the ordinary smooth case, as long as $\pi_0(t,\cdot)$ is sufficiently smooth or $K$ grows sufficiently fast, and $h_0$ decays fast enough, the sieve approximation error of our estimator $\widehat{\mu}(t)$ is asymptotically negligible and the dominating bias term is the same as that in the ordinary smooth case. The asymptotic variance is affected by the measurement error $U$. The convergence rate of the variance for bandwidth $b=h$ or $h_0$, $\{Nb\exp(-2b^{-\beta}/\gamma)\}^{-1/2}$, is degenerated by $\exp(b^{-\beta}/\gamma)$ compared to the rate $(Nb)^{-1/2}$ for the error-free case, owing to the supersmoothness of the error distribution. \label{ReplySSRate} 

From the theorem, when $h \asymp h_0$ and $\min(h,h_0) =d(\log N)^{-1/\beta}$ for a constant $d>(2/\gamma)^{1/\beta}$, one finds that the rate of variance, $N^{-1}\{e(h)\wedge e(h_0)\}^{-2}=o(h^4+h_0^4)$, is negligible compared to the asymptotic bias, and the convergence rate of $\widehat{\mu}(t)-\mu(t)$ is $\left( \log N\right)^{-2/\beta}$. This result is analogue to that in the literature on nonparametric regression with measurement error (see e.g.~\citealp{fan1991asymptotic,fan1991asymptotic,Fan1993,Delaigle2015} among others) and it achieves the optimal convergence rate for all possible nonparametric regression estimators when the regressors are measured with supersmooth errors showed in \cite{Fan1993}.

Note that under the case of supersmooth error, an explicit expression and the exact convergence rate of $\text{var}\{\eta_{h,h_0}(S,\boldsymbol{X},Y;t)\}$ is extremely hard (if not impossible) to derive without additional assumptions. In order to establish the asymptotic distribution of $\widehat{\mu}$ using Lyapunov central limit theorem, a lower bound of the deconvolution kernel's second moment is required. In particular, we require $v_h(t)\geq d_1 f_S(t)h^{d_3}\exp(2h^{-\beta}/\gamma - d_2h^{-d_4\beta})$. This is commonly imposed in the measurement error literature; see \cite{fan1991asymptotic} and \cite{Delaigle2015} among others. \cite{fan1991asymptotic} showed that this lower bound holds under some mild conditions on $\phi_U$ and $\phi_L$ (e.g.~\eqref{SuperS} and Assumption~S hold, $\phi_L(t)>c_L(1-t)^3$ for $t\in[1-\epsilon,1)$ for some $c_L,\epsilon>0$, and the real part $R_U(t)$ and the imaginary part $I_U(t)$ of $\phi_U$ satisfy $R_U(t)=o\{I_U(t)\}$ or $I_U(t)=o\{R_U(t)\}$ as $t\rightarrow\infty$). These assumptions do not exclude the usually-used kernel function $L$ defined below Assumption~S and error distributions such as Gaussian, Cauchy and Gaussian mixture.

Our explicit asymptotic linear expansion of $\widehat{\mu}(t)-\mu(t)$ in \eqref{LinearExpressionSS} is particularly helpful for statistical inference in the supersmooth error case due to the difficulty of deriving an explicit expression of the asymptotic variance. Most of the literature provides only the convergence rate; see \cite{Fan1993}, \cite{meister2006density}, and \cite{Meister2009}, among others. 

\section{Select the Smoothing Parameters}\label{sec:ParameterSelection}
In this section, we discuss how to choose the three smoothing parameters $K,h_0$, and $h$ to calculate our estimator $\widehat{\mu}(t)$ (see \eqref{pihat}, \eqref{hatG}, and \eqref{CalibrationF}). Before delving into our method, we need some preliminaries.

\subsection{Preliminaries}
The smoothing parameters in nonparametric regression are usually selected by either minimising certain cross-validation (CV) criteria or minimising an approximation of the asymptotic bias and variance of the estimator.

In nonparametric errors-in-variables regression, as pointed out by \cite{Carroll2006} and \cite{Meister2009}, approximating the asymptotic bias and variance of the estimator can be extremely challenging, if not impossible. 
Unfortunately, the CV criteria are also not computable. To see this, we assume that $K$ and $h_0$ in \eqref{pihat} and \eqref{hatG} are given for now and adapt the CV criteria to our context to choose $h$, which would be
\begin{equation}
CV(h) = \sum^N_{i=1} \big\{\widehat{\pi}(T_i,\bs{X}_i)Y_i - \widehat{\mu}^{-i}(T_i)\big\}^2 w(T_i)\,,\label{CV}
\end{equation}
where $w$ is a weight function that prevents the CV from becoming too large because of the unreliable data points from the tails of the distributions of $T$ and $\widehat{\mu}^{-i}$ denotes the estimator obtained as in \eqref{CalibrationF}, but without using the observations from individual $i$. Now, we see that \eqref{CV} is not computable in errors-in-variables regression problems since the $T_i$'s are not observable. 

To tackle this problem, \cite{Delaigle2008SIMEX} proposed combining the CV and SIMEX methods (e.g.~\citealp{Cook1994} and \citealp{Stefanski1995}). Specifically, in the simulation step, we generate two additional sets of contaminated data, namely $S_{i,d}^* = S_i + U_{i,d}^*$ and $S_{i,d}^{**} = S_{i,d}^*+U_{i,d}^{**}$, for $i=1,\ldots,N$ and $d=1,\ldots,D$ with $D$ a large number, where the $U_{i,d}^{*}$'s and $U_{i,d}^{**}$'s are i.i.d. as $U$ in \eqref{ME}. 
Now, inserting first the $S_{i,d}^*$'s and then the $S_{i,d}^{**}$'s in \eqref{pihat} and \eqref{CalibrationF} instead of the $S_i$'s, we obtain respectively $(\widehat{\pi}^*_d,\widehat{\mu}_d^*)$ and $(\widehat{\pi}^{**}_d,\widehat{\mu}_d^{**})$ for $d=1,\ldots,D$. The authors then suggested deriving two CV-type bandwidths, $\widehat{h}^*$ and $\widehat{h}^{**}$, which minimise $\sum^D_{d=1}CV_d^*(h)/D$ and $\sum^D_{d=1}CV_d^{**}(h)/D$, respectively, where
\begin{align*}
	CV_d^*(h) =& \sum^N_{i=1} \big\{\widehat{\pi}^*_d(S_i,\bs{X}_i)Y_i - \widehat{\mu}^{*,-i}_d(S_i)\big\}^2w(S_i)\,,\\
	CV_d^{**}(h) =& \sum^N_{i=1} \big\{\widehat{\pi}^{**}_d(S^*_{i,d},\bs{X}_i)Y_i - \widehat{\mu}^{**,-i}_d(S^*_{i,d})\big\}^2w(S^*_{i,d})\,,
\end{align*}
for $d=1,\ldots,D$, where $\widehat{\mu}^{*,-i}_d$ and $\widehat{\mu}^{**,-i}_d$ are obtained respectively as $\widehat{\mu}_d^*$ and $\widehat{\mu}_d^{**}$, but without using the observations from individual $i$.

The $S^{**}_{i,d}$'s are the contaminated version of the $S^*_i$'s, which is the same role as the $S^*_i$'s play to the $S_i$'s and the $S_i$'s play to the $T_i$'s. Intuitively, we then expect the relationship between $\widehat{h}^*$ and our target bandwidth $h$ to be similar to that between $\widehat{h}^{**}$ and $\widehat{h}^*$. Thus, the authors proposed an extrapolation step to obtain an estimator of $h$. Specifically, they considered that $h/\widehat{h}^\ast \approx \widehat{h}^\ast / \widehat{h}^{**}$ and used a linear back-extrapolation procedure that, in our context, would give the bandwidth
\begin{equation}
\widehat{h}_{DH} = (\widehat{h}^*)^2/\widehat{h}^{**}\,.\label{hDH}
\end{equation}

\subsection{Two-step Procedure and Local Constant Extrapolation}\label{BWSelector}
In our case, recall that we have two more smoothing parameters, $K$ and $h_0$. We can either extend the SIMEX method to choose three parameters simultaneously or choose $K$ and $h_0$ using other methods first and then apply SIMEX to choose $h$. The first option incurs a high computational burden and is unstable in practice. Thus, we adopt the second choice, which leads to a two-step procedure.

Note from Theorems~\ref{thm:normality} and \ref{thm:normality_SS} that our estimator achieves optimal rate when $h\asymp h_0$ trades off the rate of the bias $h^2+h_0^2$ and that of the standard deviation $\sqrt{v_h(t)/N+v_{h_0}(t)/N}$. Note also that the plug-in bandwidth $h_{PI}$ for the kernel deconvolution estimator with bandwidth $h$ of the density of $T$ proposed by \cite{Delaigle2002} minimises the asymptotic MSE of the estimator, whose bias is of rate of $h^2$ and standard deviation $\sqrt{v_h(t)/N}$. Thus, we should have our $h\asymp h_0 \asymp h_{PI}$. Moreover, to make $K$ satisfy all the conditions in our theorems when $h\asymp h_0$, we require $K =o(h^{-2})$.

Thus, we first set $h_0=h_{PI}$. Then, to choose $K$, we note from \eqref{ConditionalMoment} that $\mathbb{E}\{\pi_0(t,\bs{X})\exp(\bs{X})|T=t\} = \mathbb{E}\{\exp(\bs{X})\}$ holds. We propose to
choose $K = \lfloor\tilde{c} h_{PI}^{-2}\log(h_{PI}+1)\rfloor$ such that $K\geq 2$, where the constant $\tilde{c}$ minimises the following generalised CV criterion \citep{craven1978smoothing}:
\begin{equation*}
\int_{\mathcal{T}}\bigg|\frac{\sum^N_{i=1}\widehat{\pi}(t,\bs{X}_i)\exp(\bs{X}_i)L_U\{(t-S_i)/h_{PI}\}}{\sum^N_{i=1}L_U\{(t-S_i)/h_{PI}\}} - \frac{\sum^N_{i=1}\exp(\bs{X}_i)}{N}\bigg|^2 \bigg/(1-K/N)^2\,dt.%\label{Kselector}
\end{equation*} 
Such a choice of $K$ and $h_0$ is not guaranteed to minimise the error of our final estimator $\widehat{\mu}(t)$. However, with our choice of $h$ below, they guarantee the optimal convergence rate of $\widehat{\mu}(t)$ if B-spline basis is used. For the polynomial sieve basis, the smooth parameter $\alpha$ defined in Assumption~\ref{as:smooth_pi} need to be larger than 1 (see Appendix~A.5). Moreover, the simulation results showed that this works well (see Section~\ref{sec:Simulation} for more discussion).

In the second step, we could simply adopt $\widehat{h}_{DH}$ in \eqref{hDH}. However, in our numerical study, the linear back-extrapolation sometimes gave highly unstable results. We expected a larger number of $D$ to reduce the variability; for example, \cite{Delaigle2008SIMEX} used $D=20$. However, even with $D=40$, we still found some unacceptable results, which was somewhat expected, as which extrapolant function should be used in practice is unknown (\citealp{Carroll2006}, Section 5.3.2). Therefore, we introduce a new extrapolation procedure.

In particular, instead of extrapolating parametrically from $\widehat{h}^*$ and $\widehat{h}^{**}$, we suggest approximating the relationship between the $h^*_d$'s and $h^{**}_d$'s using a local constant estimator (see \citealp{Fan1996}), where $h^*_d=c^*_d h_{PI}$ and $h^{**}_d=c^{**}_d h_{PI}$ with the constants $c^*_d, c^{**}_d$ minimise $CV_d^*(h)$ and $CV_d^{**}(h)$, respectively, for $d=1,\ldots,D$. Then, we take this approximated relationship as the extrapolant function. Specifically, we choose the bandwidth $h$ to be
\begin{equation}
	\widehat{h} = \frac{\sum^{D}_{d=1}h_d^*\cdot\varphi\{(\widehat{h}^* - h^{**}_d)/b\}}{\sum^D_{d=1}\varphi\{(\widehat{h}^* - h^{**}_d)/b\}}\,,\label{hhat}
\end{equation}
where $\varphi$ is the Gaussian kernel function. The bandwidth $b$ here is selected by leave-one-out cross-validation. Local constant estimator has been well studied and
widely used, and can work fairly fast and stable. In our simulation study, we found that $D=35$ is sufficiently large to ensure good performance.

\begin{rem}\label{rem:TruncateLU}
	Recall from \eqref{hatG} that some of the deconvolution kernel $L_U\{(t-S_i)/h_0\}$'s may take negative values, making the maximisation of $\widehat{G}_t(\lambda)$ not strictly concave in finite samples. With $h_0=h_{PI}$, truncating those negative $L_U\{(t-S_i)/h_0\}$'s to 0 is a fast and stable way to solve the problem. The simulation performed well.
\end{rem}

\section{Numerical Properties}\label{sec:Simulation}
\begin{figure}[t]
	\begin{center}
	\includegraphics[width=.32\textwidth]{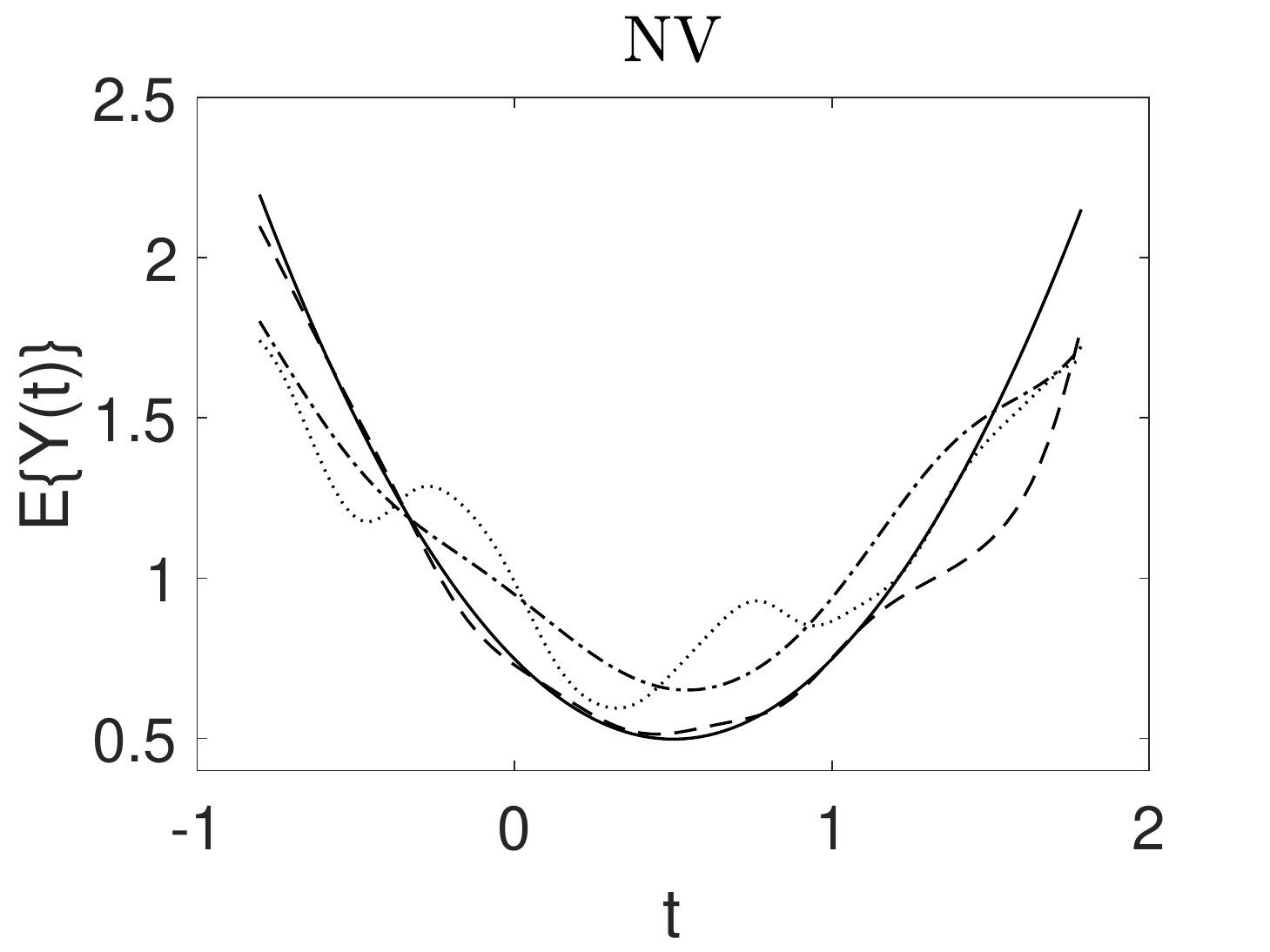}
	\includegraphics[width=.32\textwidth]{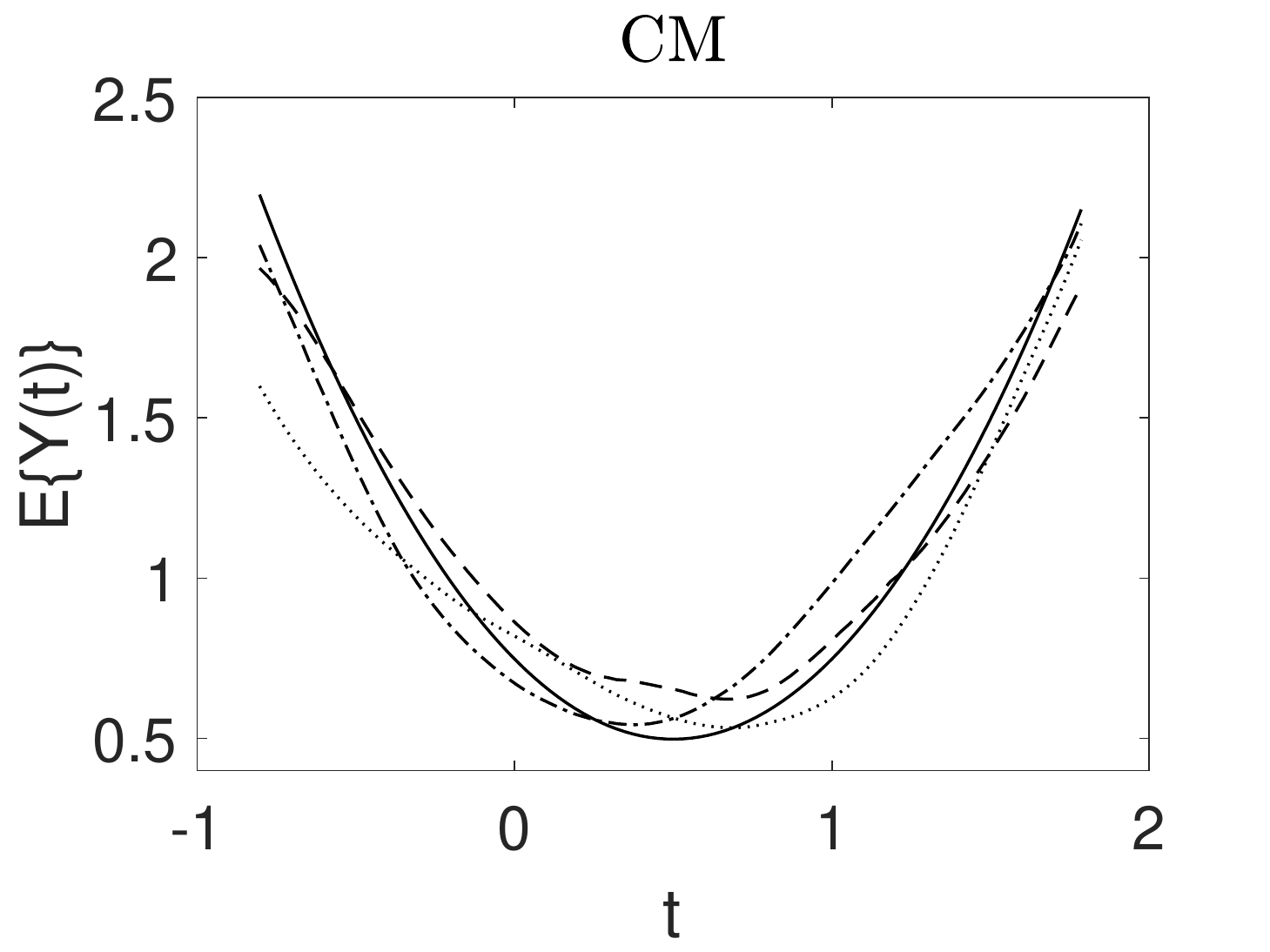}
	\includegraphics[width=.32\textwidth, height=3cm]{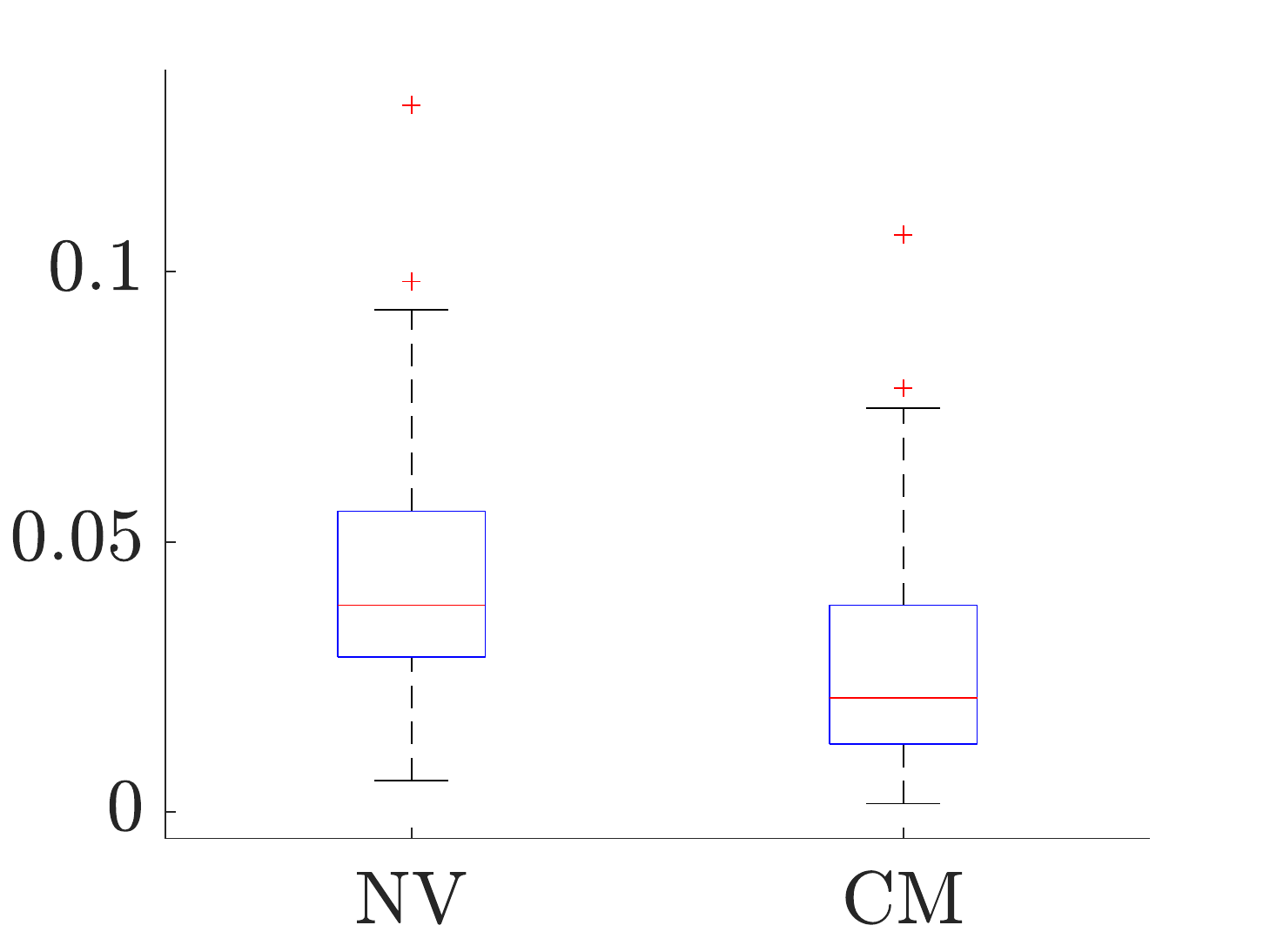}

	\vspace{.1cm}

	\includegraphics[width=.32\textwidth]{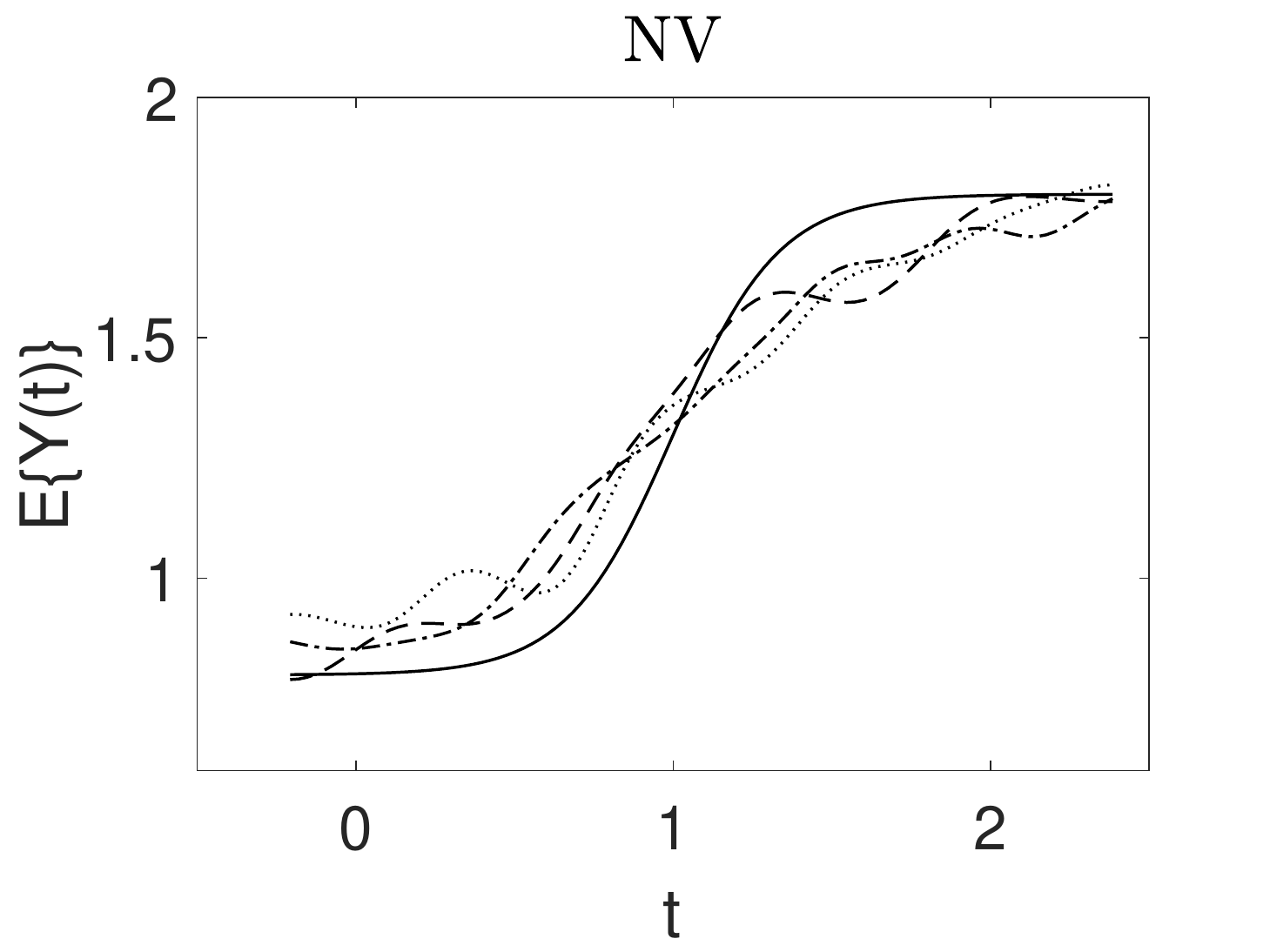}
	\includegraphics[width=.32\textwidth]{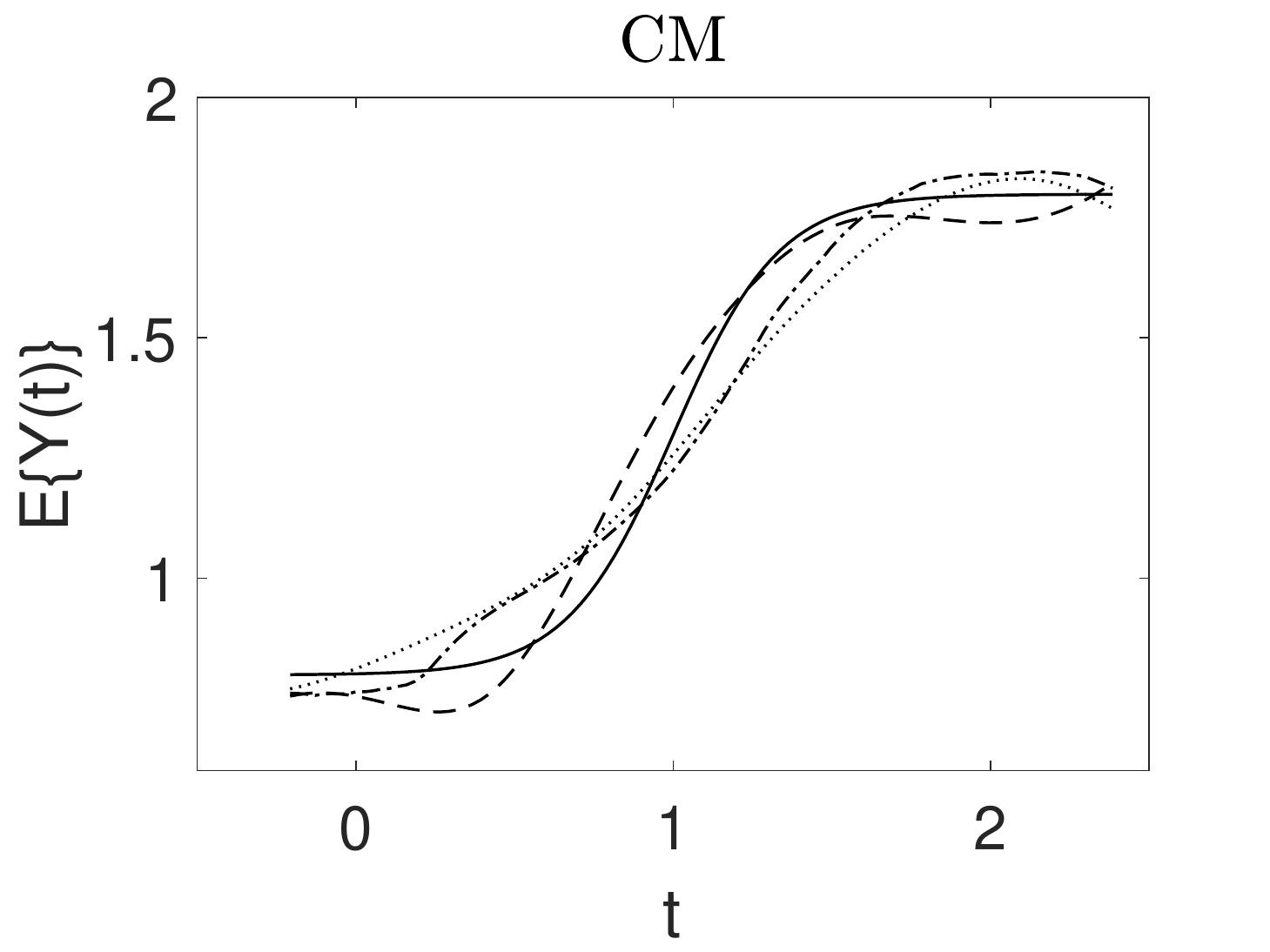}
	\includegraphics[width=.32\textwidth, height=3cm]{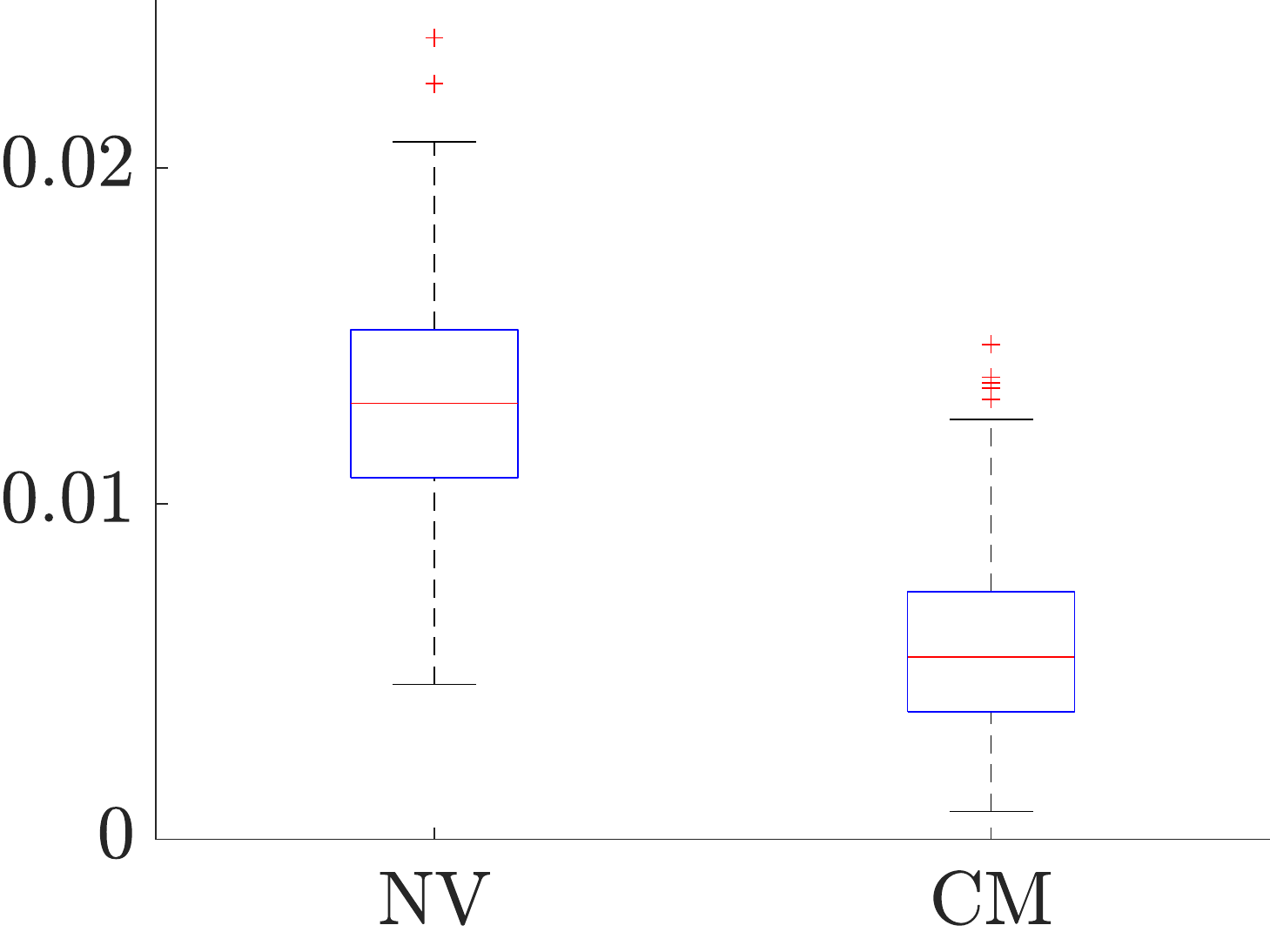}
	
	\end{center}
\caption{Plots of the true curve (solid line) and estimated curves corresponding to the 1st (dashed line), 2nd (dash-dotted line), and 3rd (dotted line) quartiles of the ISEs from the 200 Monte Carlo samples of models~1 (row 1) and 2 (row 2) with Laplace measurement errors and $N=500$. The third column depicts the boxplots of the 200 ISEs of each estimator.} \label{Fig:1}
\end{figure}

\begin{figure}[t]
	\begin{center}
	\includegraphics[width=.32\textwidth]{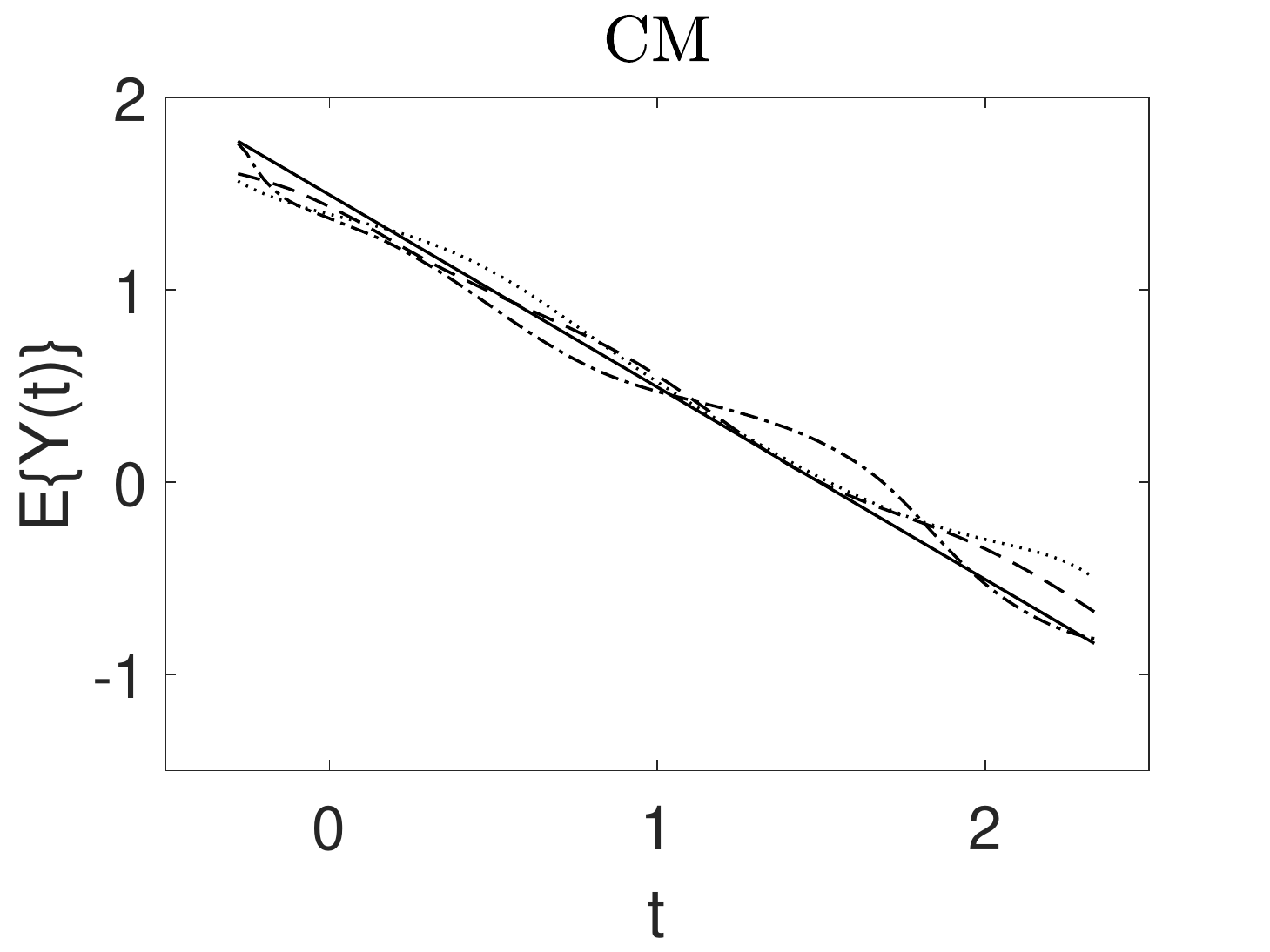}
	\includegraphics[width=.32\textwidth]{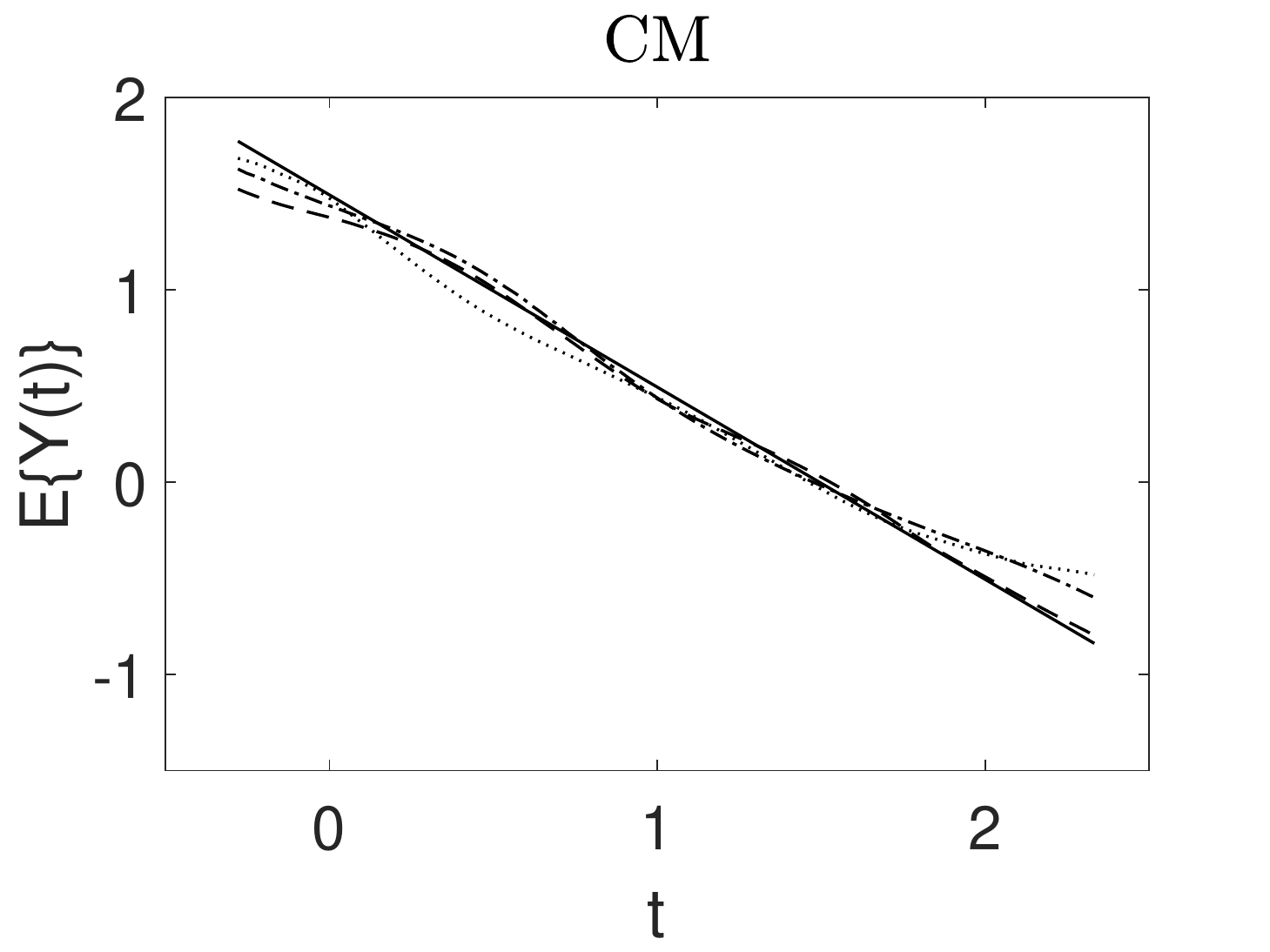}
	\includegraphics[width=.32\textwidth, height=3cm]{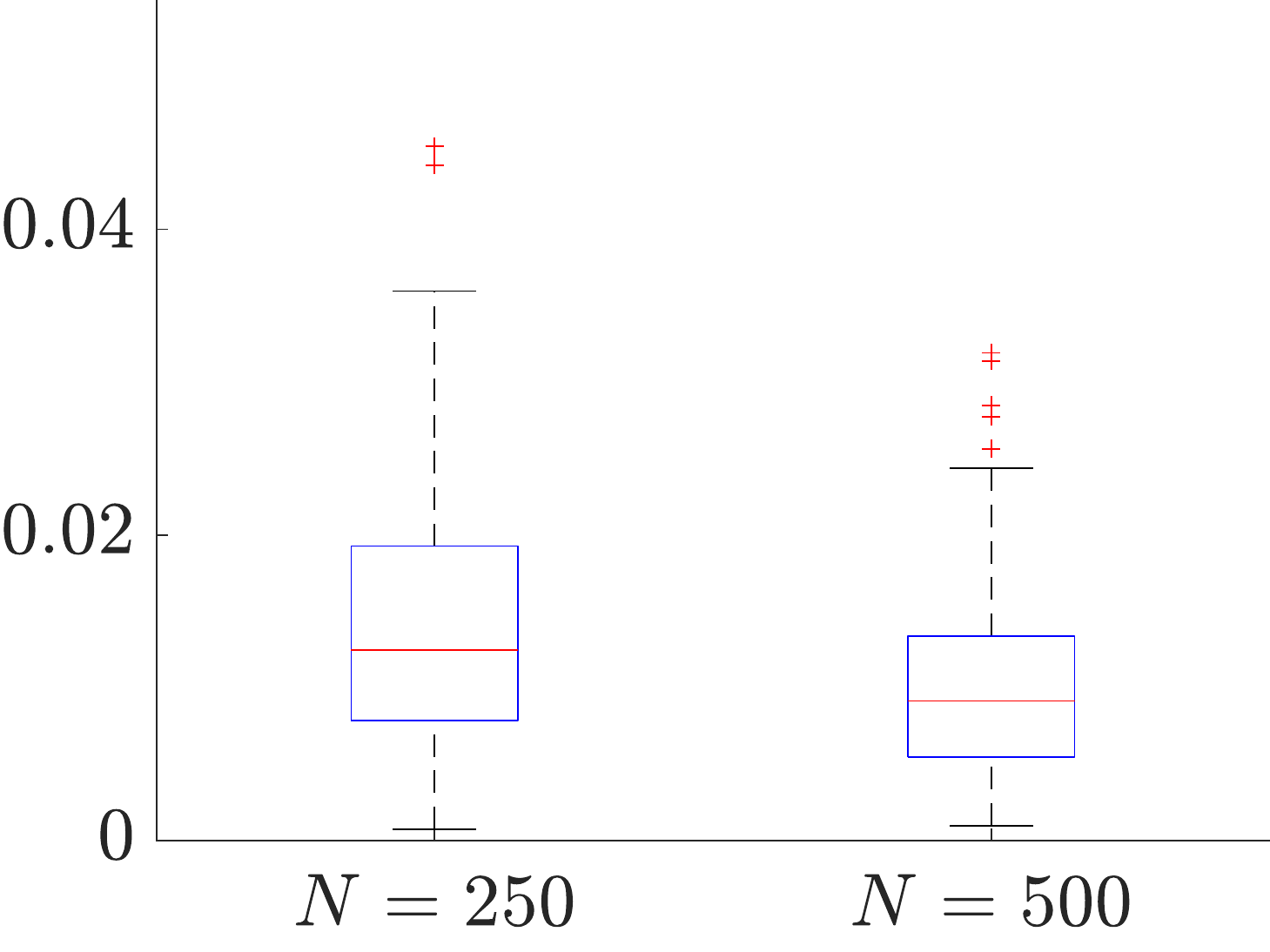}

	\vspace{.1cm}

	\includegraphics[width=.32\textwidth]{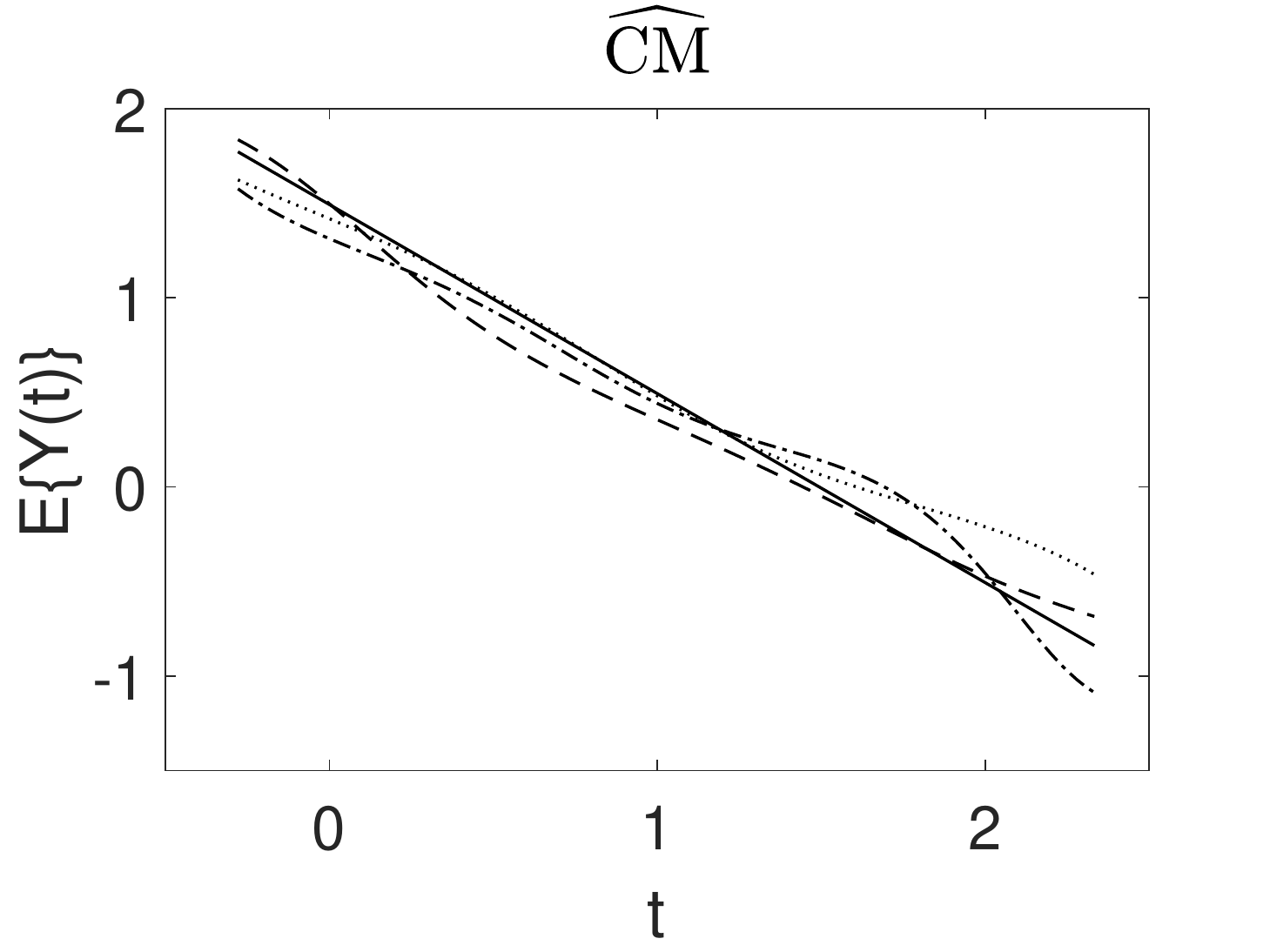}
	\includegraphics[width=.32\textwidth]{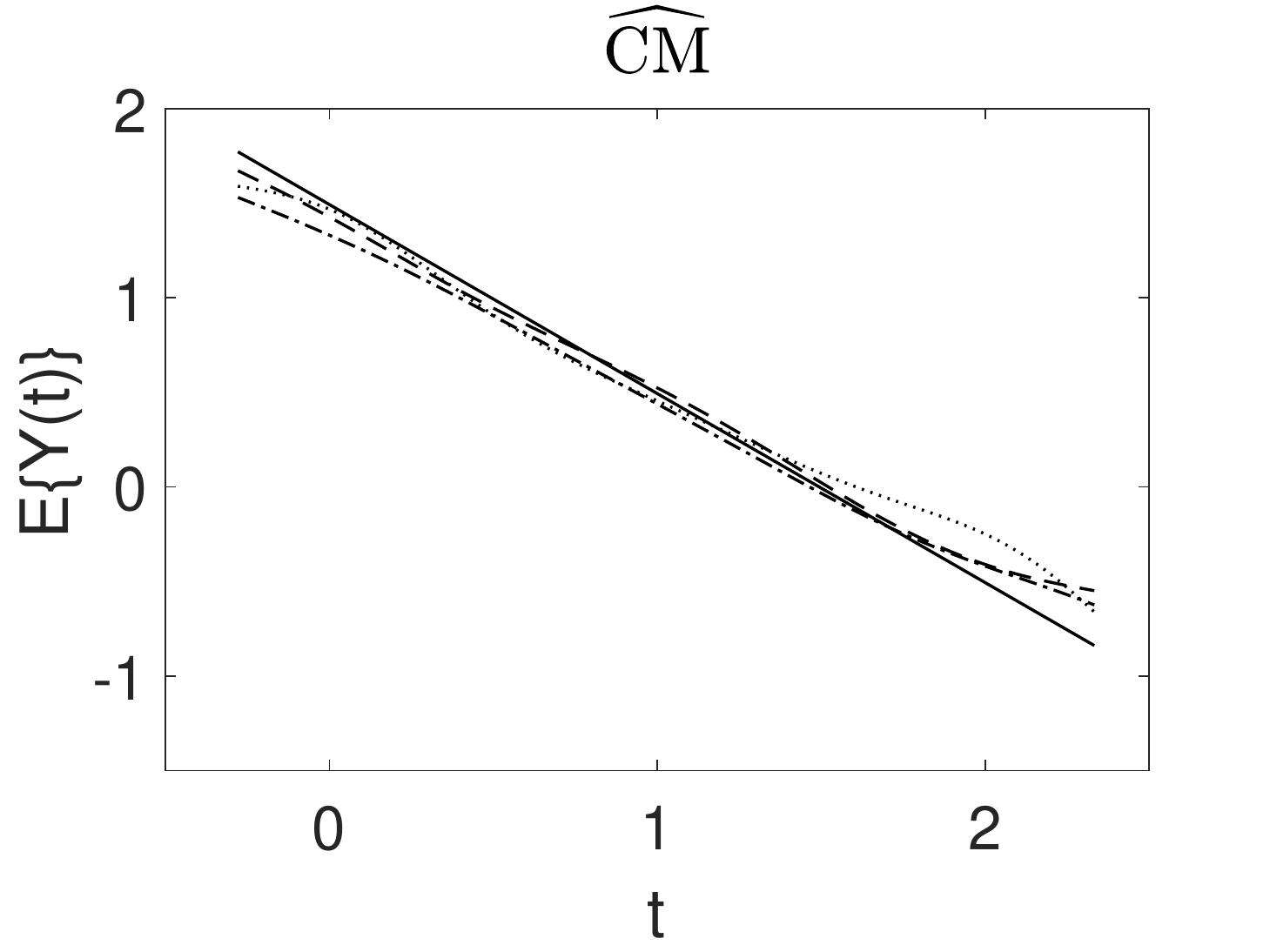}
	\includegraphics[width=.32\textwidth, height = 3cm]{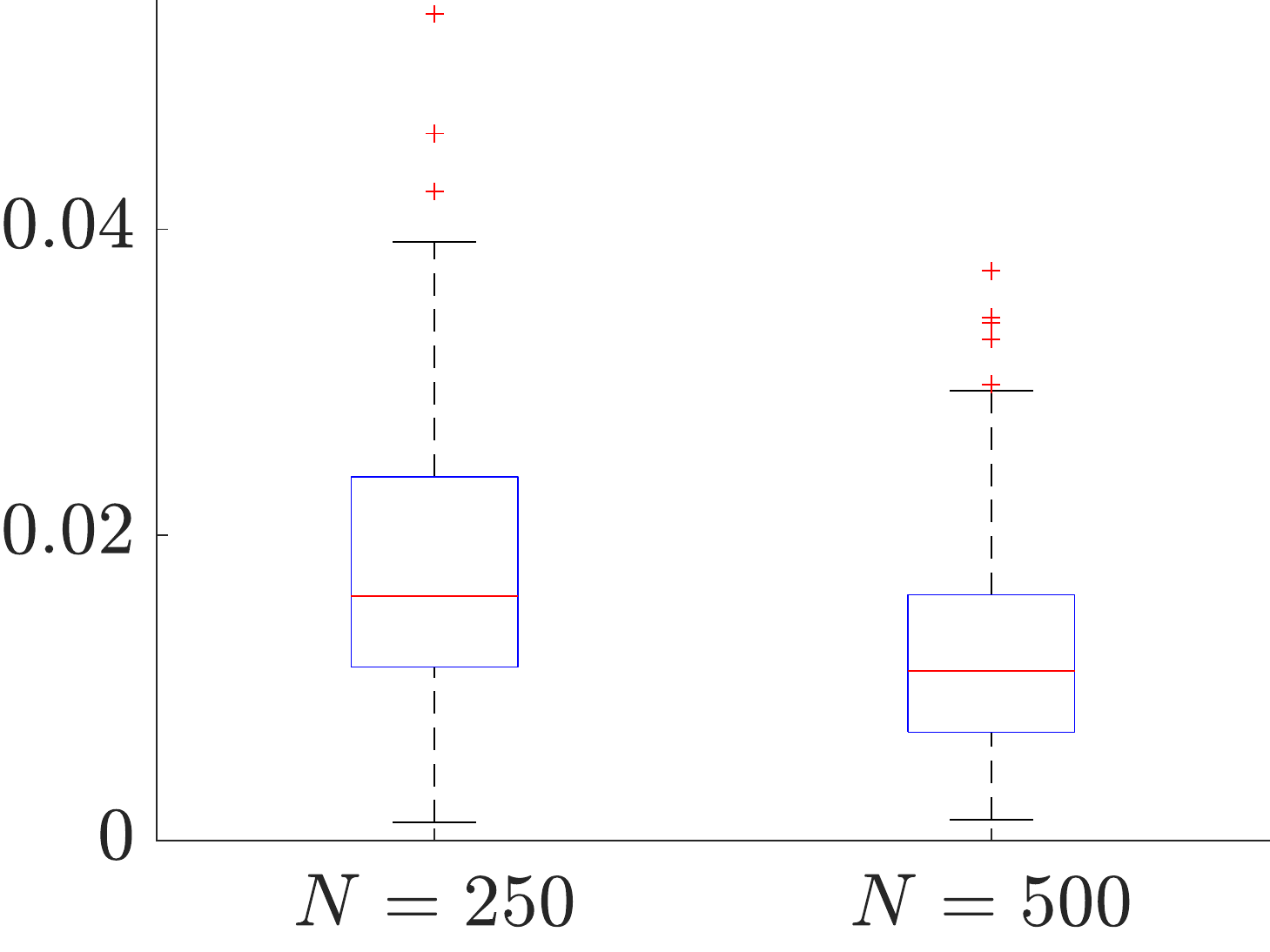}
	\end{center}
\caption{Plots of the true curve (solid line) and estimated curves CM (row 1) and $\widehat{\text{CM}}$ (row 2) corresponding to the 1st (dashed line), 2nd (dash-dotted line), and 3rd (dotted line) quartiles of the 200 ISEs from model~3 with Laplace measurement errors and $N=250$ (left) and $N=500$ (centre) and the boxplots of the 200 ISEs (right).}\label{Fig:2}
\end{figure}

\begin{figure}[t]
	\begin{center}
	\includegraphics[width=.32\textwidth]{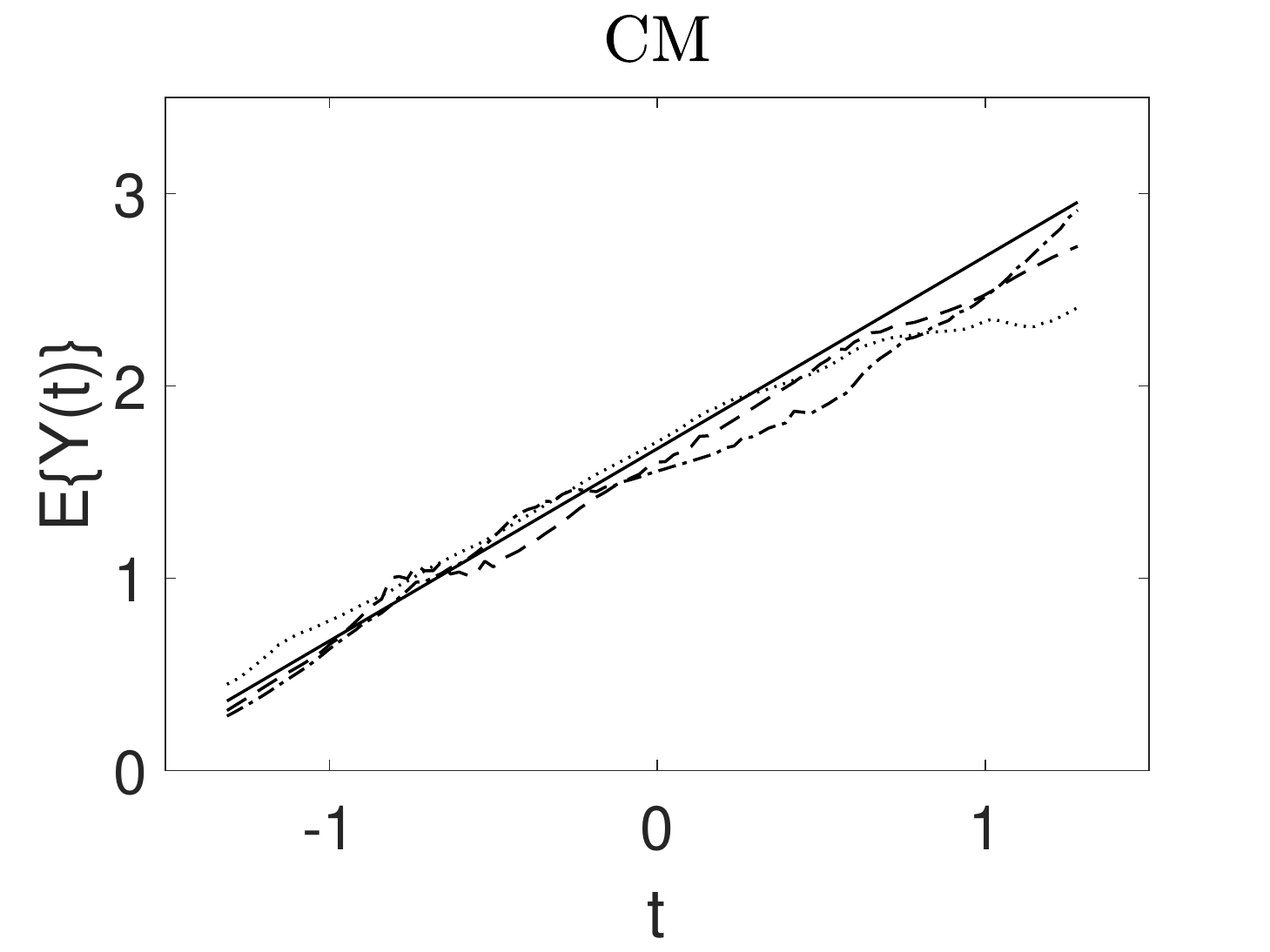}
	%\hspace*{-.4cm}
	\includegraphics[width=.32\textwidth]{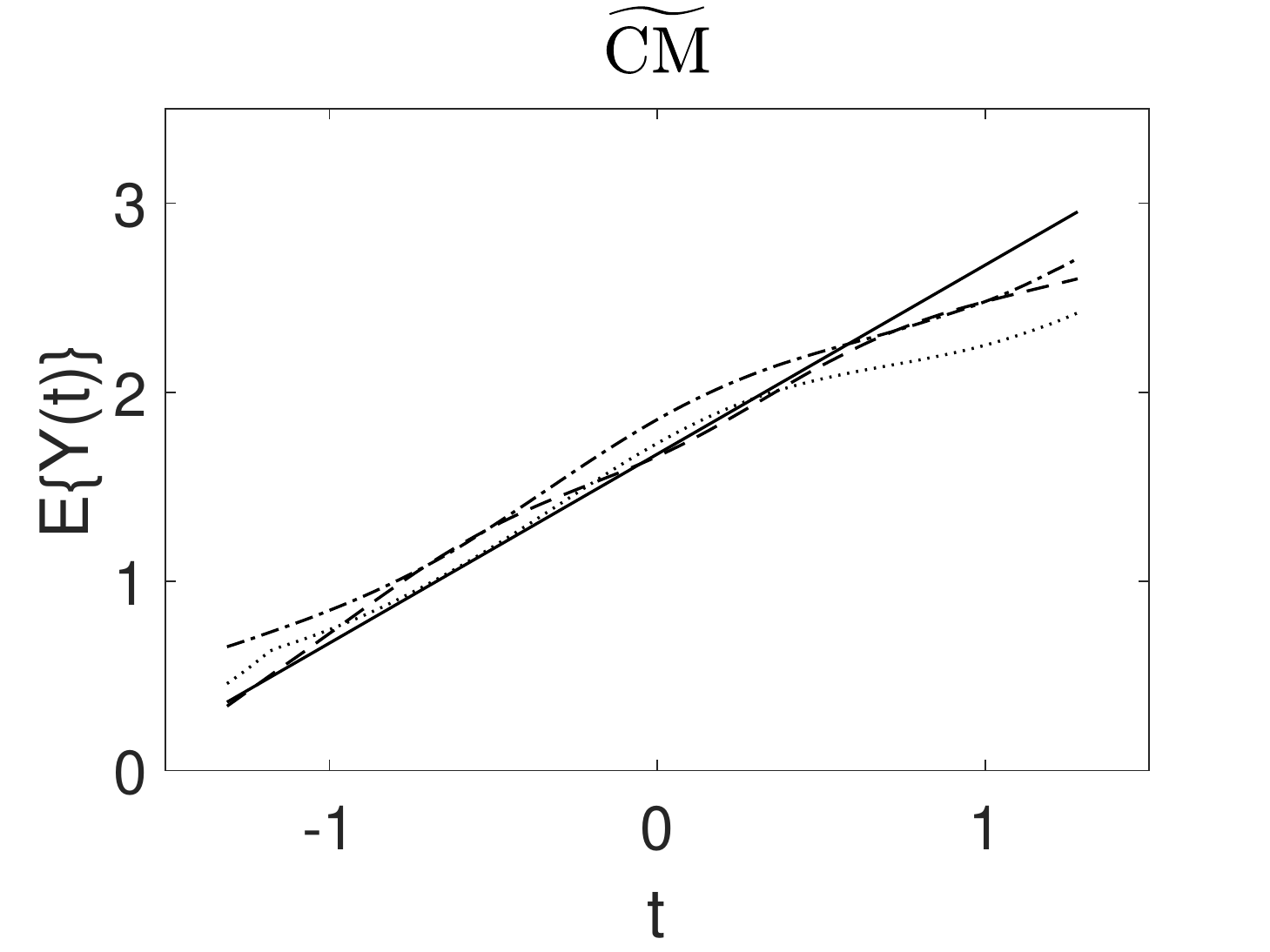}
	%\hspace*{-.4cm}
	\includegraphics[width=.32\textwidth]{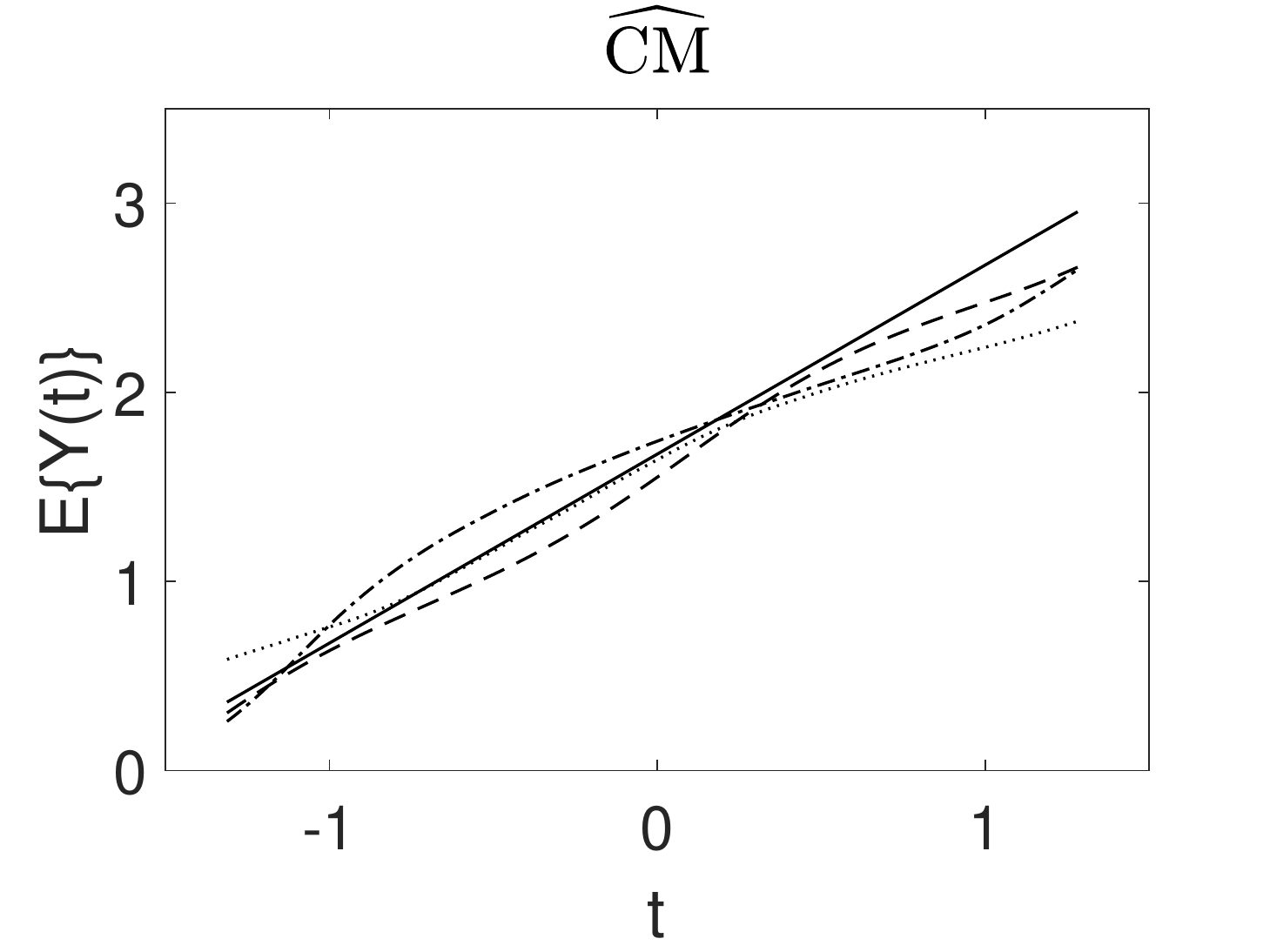}\end{center}
\caption{Plots of the true curve (solid line) and estimated curves corresponding to the 1st (dashed line), 2nd (dash-dotted line), and 3rd (dotted line) quartiles of the 200 ISEs from model~4 with Gaussian measurement errors and $N=250$.} \label{Fig:3}
\end{figure}

\subsection{Simulation Settings}
Let $\xi_{x,1},\ldots,\xi_{x,10}$ be i.i.d. uniform random variables supported on $[0,1]$ and $\xi_t\sim N(0,1)$. We consider the following four models in which $\xi_y$ is generated from a standard normal distribution for models 1 and 4 and from a uniform distribution supported on $[0,1]$ for models 2 and 3:
\begin{enumerate}[1.]
	\item $X = 0.3+0.4\xi_{x,1}$, $T=X+\xi_t$, and $Y^\ast(t)=(t-0.5)^2+X+\xi_y$ ($X$ affects $T$ and $Y$ linearly);
	\item $X = \sum^2_{j=1}0.3\xi_{x,j}$, $T=1+X^2+\xi_t$, and $Y^\ast(t)=\exp(-6+6t)/\{1+\exp(-6+6t)\}+X+\xi_y$ ($X$ affects $T$ nonlinearly and $Y$ linearly);
	\item $X = \sum^{10}_{j=1}0.2\xi_{x,j}$, $T=X+\xi_t$, and $Y^\ast(t)=-t+\sqrt{X}+\xi_y$ ($X$ affects $T$ linearly but $Y$ nonlinearly);
	\item $X = 0.2+0.6\xi_{x,1}$, $T=\sqrt{X}-0.7+\xi_t$, and $Y^\ast(t)=t+\exp(X)+\xi_y$ ($X$ affects $T$ and $Y$ nonlinearly).
\end{enumerate}
For each model, we generate 200 samples of $(S,X,Y)$ of size 250 or 500, where $S=T+U$ with $U$ either a Laplace random variable with mean 0 and $\text{var}(U)/\text{var}(T)=0.25$ or a mean zero Gaussian random variable with $\text{var}(U)/\text{var}(T)=0.2$.

For each combination of the model, sample size, and measurement error type, we calculate our estimator in \eqref{CalibrationF}. To measure the quality of the estimator, we calculate the integrated squared errors $\text{ISE}=\int_{q_{0.1}}^{q_{0.9}} \big\{\widehat{\mu}(t) - \mu(t)\big\}^2\,dt$, where $q_{0.1}, q_{0.9}$ are the 10th and 90th quantiles of $T$, respectively. 

To highlight the importance of considering the measurement errors in the estimation, we also calculate the naive estimator that ignored the error for each sample. That is, we apply the estimator of \citet{Ai_Linton_Motegi_Zhang_cts_treat, Ai2020} to our data by replacing the $T_i$'s there with the $S_i$'s. Specifically, the naive estimator is  
\begin{equation}
\frac{\sum^N_{i=1}\widetilde{\pi}(S_i,\bs{X}_i)Y_i\varphi\{(t-S_i)/h_n\}}{\sum^N_{i=1}\varphi\{(t-S_i)/h_n\}}\,, \quad t\in\mathcal{T}\,,\label{NaiveDef}
\end{equation}
where $\varphi$ is the standard normal density function and $\widetilde{\pi}(S_i,\bs{X}_i)$ is calculated using \eqref{AipiEst} (see \citealp{Ai_Linton_Motegi_Zhang_cts_treat, Ai2020} for more details).

Unless otherwise specified, we take the kernel function $L$ for the deconvolution kernel method to be the one whose Fourier transform is given by 
$
\phi_{L}(u) = (1-u^2)^3\cdot \mathbbm{1}_{[-1,1]}(u).
$

%= \rho^{'}\{u_{K_1}^\top(S_i)\widehat{\Gamma}_{K_1\times K_2}v_{K_2}(\bs{X}_i)\}$ for $i=1,\ldots,N$, $\rho$ defined below \eqref{hatG} and $\widehat{\Gamma}_{K_1\times K_2}$ maximises 
%$$
% \frac{1}{N}\sum^N_{i=1} \rho\big\{u_{K_1}^\top(S_i)\widehat{\Gamma}_{K_1\times K_2}v_{K_2}(\bs{X}_i)\big\} - \bigg\{\frac{1}{N}\sum^N_{i=1}u_{K_1}(S_i)\bigg\}^\top\Gamma \bigg\{\frac{1}{N}\sum^N_{i=1}v_{K_2}(\bs{X}_i)\bigg\}\,,
%$$
%with $K_1,K_2$ the number of basis functions to be chosen.

%We also calculated the plug-in estimator (PI) defined in \eqref{PlugInF} and the outcome regression estimator (OR) defined in \eqref{ORModel} for comparison.

%Note that both the PI and OR methods have more than one smoothing parameters to choose. Applying the SIMEX method to choose the parameters simultaneously is again computationally expensive. We suggest to choose $h_1=h_2=h_{PI}$ defined in section~\ref{BWSelector} for the PI method, and $\bs{H}$ to be the bandwidth proposed by \cite{Bowman1997} for both PI method and OR method, and then apply the modified SIMEX method in section~\ref{BWSelector} to choose $h$ in \eqref{PlugInF} and in \eqref{mtX}.

To illustrate the potential benefit of using our methods over the naive estimator without confounding the effect of the smoothing parameter selectors, we first use the theoretically optimal smoothing parameters for each method. These parameters simultaneously minimise the integrated squared error (ISE) for each method, resulting in the optimal naive estimator (NV) and the proposed conditional moment estimator $\widehat{\mu}$ (CM).

Recall from Section~\ref{BWSelector} that we do not choose $K$ and $h_0$ by minimising the estimated ISE of our estimator. To see how much we might lose by doing so, we calculate the estimator with our choice of $K$ and $h_0$ and the optimal $h$ that minimised the ISE, which is denoted by $\widetilde{\text{CM}}$.

Finally, to assess the performance of our method in practice, we calculate $\widehat{\mu}$ using the smoothing parameters selected from the data using the method in Section~\ref{BWSelector} and denote it by $\widehat{\text{CM}}$. We take the weight function $w$ to be an indicator function that equals 1 when the $S_i$'s or $S^*_{i,d}$'s are within their 5\% to 95\% quantiles, and 0 otherwise. We also compute the naive estimator in \eqref{NaiveDef} with the $K_1, K_2$, and $h_n$ selected using the 10-fold CV method, denoted by $\widehat{\text{NV}}$, to make a comparison.

\subsection{Simulation Results}

In this section, we show our simulation results. The full simulation results of the 200 values of the ISE of each estimator obtained from the 200 simulated samples for models~1 to 4 can be seen in the boxplots in Appendix~A.6. Figures~\ref{Fig:1} to \ref{Fig:3} depict the true curve of the model and three estimated curves corresponding to the 1st, 2nd, and 3rd quartiles of the 200 ISE values.

Overall, the simulation results show that our methods with the theoretically optimal smoothing parameters (CM) perform better than that of the naive one (NV). This confirms the advantages of our methods over the naive one by adapting the estimation to the measurement errors. A graphical example is presented in Figure~\ref{Fig:1}, which shows the quartile curves of NV and CM for models~1 and 2 with Laplace measurement errors and $N=500$.

The simulation results also confirm our theoretical proportion that the performance of our method improves as the sample size increases. Figure~\ref{Fig:2} exemplifies the effect of increasing $N$ by depicting the quartile curves and ISE boxplots of CM and $\widehat{\text{CM}}$ for model~3 when the measurement errors follow a Laplace distribution with $N=250$ and $N=500$. The improvement with the increase in sample size can also be seen in the boxplots in the Appendex~A.6.\label{ReplySampleSize} %For comparison, we also show the quartile curves of NV estimator there. As a biased estimator, the improvement of NV estimator with increasing sample size is marginal, especially near the boundary.

Comparing the ISE values of $\widetilde{\text{CM}}$ with those of CM, we find that our choice of $K$ and $h_0$ discussed in Section~\ref{BWSelector} lowers the performance of our estimator only marginally in most cases. Recall that $K$ is the number of polynomial basis functions of $\bs{X}$ used to estimate $\pi_0$ and $h_0$ is the bandwidth used to estimate $\mathbb{E}\{\pi_0(t,\bs{X})|T=t\}$, which are related to the relationship between $\bs{X}$ and $T$ as well as that between $\bs{X}$ and $Y^*(t)$. We thus consider the nonlinear relationship between $T$ and $\bs{X}$ or between $Y^*(t)$ and $\bs{X}$ in the simulation models (see~models~2 to 4). Our choice of $K$ and $h_0$ still works well compared with the optimal one. Figure~\ref{Fig:3} provides a graphical example from model~4, where $\bs{X}$ affects both $T$ and $Y^*(t)$ nonlinearly, with Gaussian measurement errors and $N=250$.

Finally, from the behaviour of $\widehat{\text{CM}}$ in all our simulation studies, we observe that our modified SIMEX method introduced in Section~\ref{BWSelector} performs well and stably. This can also be seen from Figures~\ref{Fig:2} and \ref{Fig:3}.

\subsection{Real Data Example}\label{sec:RealData}
We demonstrate the practical value of our data-driven $\widehat{\text{CM}}$ estimator using Epidemiologic Study Cohort data from NHANES-I. We estimate the causal effect of the long-term log-transformed daily saturated fat intake on the risk of breast cancer based on a sample of 3,145 women aged 25 to 50. The data were analysed by \cite{Carroll2006} using a logistic regression calibration method, and they are available from \url{https://carroll.stat.tamu.edu/data-and-documentation}. The daily saturated fat intake was measured using a single 24-hour recall. Specifically, the log-transformation was taken as $\log(5+\text{saturated fat})$. Previous nutrition studies have estimated that over 75\% of the variance in those data is made up of measurement error. According to \cite{Carroll2006}, it is reasonable to assume the classical measurement error model (i.e. \eqref{ME}), with a Gaussian measurement error $U$ on the data. The outcome variable $Y$ takes 1 if the individual has breast cancer and 0 otherwise. The covariates in $\bs{X}$ are age, the poverty index ratio, the body mass index, alcohol use (yes or no), family history of breast cancer, age at menarche (a dummy variable taking 1 if age is $\leq 12$), menopausal status (pre or post), and race, which are assumed to have been measured without appreciable error.

We first apply our estimator to the data for a Gaussian measurement error with an error variance $\text{var}(U)/\text{var}(S)=0.75$. That corresponds to $\text{var}(U)/\text{var}(T) = 3$. As pointed out by \cite{delaigle2004bootstrap}, the error variances estimated by other nutrition studies may be inaccurate. Thus, we also consider cases in which $\text{var}(U)/\text{var}(S) = 0.43$, $\text{var}(U)/\text{var}(S) = 0.17$, and $\text{var}(U)=0$ (i.e. $\text{var}(U)/\text{var}(T) = 0.75$, $\text{var}(U)/\text{var}(T) = 0.2$ and the error-free case).

\begin{figure}[t]
\centering
\includegraphics[width = .4\textwidth]{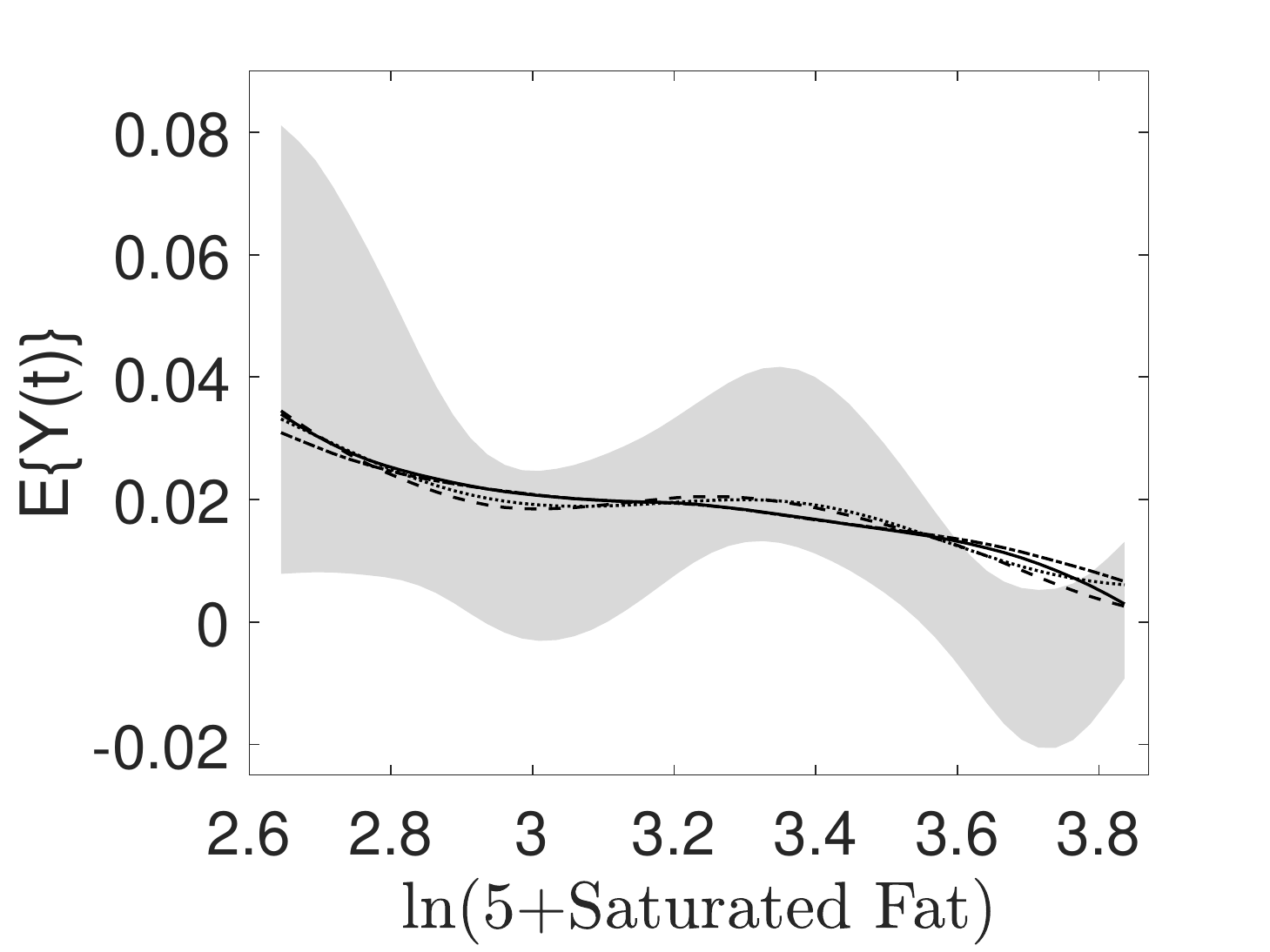}
\caption{Estimation of the treatment effect of the log-saturated fat intake on the risk of breast cancer for a Gaussian error of $\text{var}(U)/\text{var}(S) =0.75$ (solid line), $\text{var}(U)/\text{var}(S) = 0.43$ (dash-dotted line), $\text{var}(U)/\text{var}(S) = 0.17$ (dashed line) with its 95\% pointwise confidence band (shaded), and $\text{var}(U)=0$ (dotted line).}\label{Fig:NHANES}
\end{figure}

Figure~\ref{Fig:NHANES} presents the estimated curves using the smoothing parameters selected as described in section~\ref{sec:ParameterSelection} and a 95\% undersmoothing pointwise confidence band for $\text{var}(U)/\text{var}(S)=0.17$ (see Appendix~A.4 for the method and the confidence bands for $\text{var}(U)/\text{var}(S)=0.43$ and 0.75). Overall, the estimated risk of breast cancer shows a decreasing trend across the range of transformed saturated fat intake. When the measurement error variance is 0.17 of $\text{var}(S)$ or 0, there is a marginal increasing trend between $t=3$ and 3.4. The 95\% confidence bands for $\text{var}(U)/\text{var}(S)=0.17$ and 0.43 show an overall decreasing trend with a slight increase between $t=3$ and 3.4.  These findings concur with the results of \cite{Carroll2006}, who found in their multivariate logistic regression calibration that the coefficient of the log-transformed saturated fat intake on the risk of breast cancer was significant and negative. However, the results should be treated with extreme caution because of possible misclassification in the breast cancer data and the lack of follow-up of breast cancer cases with high fat intakes; see \citet[Chap 3.3]{Carroll2006}. 

%Indeed, it is impossible to conclude whether the saturated fat intake has a significant effect on the risk of breast cancer or if the risk decreases with the intake monotonically based on these results; hence, it would be interesting for future studies to address these problems.

%\section{Conclusion Remarks}
%{\color{red}Refer to \cite{Delaigle2008} on unknown measurement error density and heteroscedastic measurement errors (Or put them where we made known measurement error distribution assumption).}

\section*{Acknowledgements}
The authors would like to sincerely thank the  Steffen Lauritzen, the Associate Editor, and the two
referees for their constructive suggestions and comments. Wei Huang's research was supported by the Professor Maurice H. Belz Fund of the University of Melbourne.
Zheng Zhang is
supported by the fund from the National Natural Science Foundation of China
[grant number 12001535], Natural Science Foundation of Beijing [grant number
1222007], and the fund for building world-class universities (disciplines)
of Renmin University of China [project number KYGJC2022014]. The authors contributed equally to this work and are listed in the alphabetical order.

\bibliographystyle{rss}
\bibliography{reference}

\end{document}